\documentclass[]{article}

\usepackage{todonotes}

\usepackage{graphicx}
\usepackage{amsfonts}
\usepackage{amsmath}
\usepackage{amssymb}
\usepackage{fancyhdr}
\usepackage{titlesec}
\usepackage{indentfirst}
\usepackage{booktabs}
\usepackage{verbatim}
\usepackage{color}
\usepackage{amsthm}
\usepackage{url}
\usepackage{cite}

\usepackage[page,header]{appendix}
\usepackage{titletoc}

\numberwithin{equation}{section}
\newtheorem{theorem}{Theorem}[section]

\newtheorem{remark}[theorem]{Remark}

\newcommand{\rd}{\mathrm{d}}

\begin{document}

\title{A fast spectral method for the Boltzmann collision operator with general collision kernels%
	\footnote{Los Alamos Report LA-UR-16-26555.Funded by the Department of Energy at Los Alamos National Laboratory under contract DE-AC52-06NA25396. }
		\footnote{This manuscript has been authored, in part, by UT-Battelle, LLC, under Contract No. DE-AC0500OR22725 with the U.S. Department of Energy. The United States Government retains and the publisher, by accepting the article for publication, acknowledges that the United States Government retains a non-exclusive, paid-up, irrevocable, world-wide license to publish or reproduce the published form of this manuscript, or allow others to do so, for the United States Government purposes. The Department of Energy will provide public access to these results of federally sponsored research in accordance with the DOE Public Access Plan (\texttt{http://energy.gov/downloads/doe-public-access-plan}).}
}

\author{Irene M. Gamba\footnote{Department of Mathematics and Institute for Computational Engineering and Sciences, The University of Texas at Austin, Austin, TX 78712, USA (gamba@math.utexas.edu). I. Gamba's research was partially supported  by NSF grant DMS-1413064  and NSF RNMS (KI-Net) grant DMS-1107465.},
	Jeffrey R. Haack\footnote{Computational Physics and Methods Group, Los Alamos National Laboratory, Los Alamos, NM 87545, USA (haack@lanl.gov). J. Haack's research was partially supported by NSF grant DMS-1109625 and NSF RNMS (KI-Net) grant DMS-1107465.}, 
	Cory D. Hauck\footnote{Computational and Applied Mathematics Group, Oak Ridge National Laboratory, Oak Ridge, TN 37831, USA (hauckc@ornl.gov). C. Hauck's research was supported by the U.S. Department of Energy, Office of Science, Office of Advanced Scientific Research.}, 
	and 
	Jingwei Hu\footnote{Department of Mathematics, Purdue University, West Lafayette, IN 47907, USA (jingweihu@purdue.edu). J. Hu's research was partially supported by NSF grant DMS-1620250 and a startup grant from Purdue University. }
	} 
\maketitle


\begin{abstract}
We propose a simple fast spectral method for the Boltzmann collision operator with general collision kernels. In contrast to the direct spectral method \cite{PR00, GT09} which requires $O(N^6)$ memory to store precomputed weights and has $O(N^6)$ numerical complexity, the new method has complexity $O(MN^4\log N)$, where $N$ is the number of discretization points in each of the three velocity dimensions and $M$ is the total number of discretization points on the sphere and $M\ll N^2$. Furthermore, it requires no precomputation for the variable hard sphere (VHS) model and only $O(MN^4)$ memory to store precomputed functions for more general collision kernels. Although a faster spectral method is available \cite{MP06} (with complexity $O(MN^3\log N)$), it works only for hard sphere molecules, thus limiting its use for practical problems. Our new method, on the other hand, can apply to arbitrary collision kernels. A series of numerical tests is performed to illustrate the efficiency and accuracy of the proposed method.
\end{abstract}

{\small 
{\bf Key words.} Boltzmann collision integral, spectral method, convolution, fast Fourier transform, Lebedev quadrature.

{\bf AMS subject classifications.} 35Q20, 65M70.

}



\section{Introduction}

Kinetic theory describes the non-equilibrium dynamics of a gas or any system comprised of a large number of particles. When well-known fluid mechanical laws of Navier-Stokes and Fourier become inadequate, kinetic equations provide rich information at the mesoscopic level and have found applications in various fields such as rarefied gas dynamics \cite{Cercignani00}, radiative transfer \cite{Chandrasekhar}, semiconductor modeling \cite{MRS}, and biological and social sciences \cite{NPT}. Our main focus in this paper is the Boltzmann equation which constitutes the central model in kinetic theory and takes the form \cite{Cercignani, CC, Villani02}:
\begin{equation} \label{CBE}
\frac{\partial f}{\partial t}+v\cdot \nabla_x f=\mathcal{Q}(f), \quad t>0, \ x\in \Omega\subset\mathbb{R}^3, \ v\in \mathbb{R}^3.
\end{equation}
Here $f = f(t,x,v)$ is the phase space distribution function, which depends on time $t$, position $x$, and particle velocity $v$; and $\mathcal{Q}$ is the Boltzmann collision operator, which models binary interactions between particles:%
\footnote{The variables $t$ and $x$ are suppressed because $\mathcal{Q}$ acts on $f$ only through the velocity.}
\begin{equation} \label{CO}
\mathcal{Q}(f)(v)=\int_{\mathbb{R}^3}\int_{S^2}\mathcal{B}(v-v_*,\omega)\left[f(v')f(v_*')-f(v)f(v_*)\right]\,\rd{\omega}\,\rd{v_*}.
\end{equation}
In this formula, $(v',v_*')$ and $(v,v_*)$ represent the velocity pairs before and after a collision.  The requirement that momentum and energy are conserved during such a collision means that $(v',v_*')$ can be expressed in terms of $(v,v_*)$:
\begin{equation} 
v'=\frac{v+v_*}{2}+\frac{|v-v_*|}{2}\omega, \quad  v_*'=\frac{v+v_*}{2}-\frac{|v-v_*|}{2}\omega,
\end{equation}
where the parameter $\omega$ varies over the unit sphere $S^2$. The collision kernel $\mathcal{B}$ is a non-negative function that depends on its arguments only through $|v-v_*|$ and cosine of the deviation angle $\theta$ (the angle between $v-v_*$ and $v'-v_*'$). Thus $\mathcal{B}$ is often written as
\begin{equation} \label{CK}
\mathcal{B}(v-v_*, \omega)=B(|v-v_*|,\cos \theta), \quad \cos\theta=\frac{\omega\cdot (v-v_*)}{|v-v_*|}.
\end{equation}
The specific form of $B$ can be determined from the intermolecular potential using scattering theory \cite{Cercignani}. For numerical purposes, a commonly used collision kernel is the variable hard sphere (VHS) model proposed by Bird \cite{Bird}:
\begin{equation} \label{VHS}
B=b_{\gamma}|v-v_*|^{\gamma}, \qquad b_{\gamma} >0, \quad  0\leq\gamma\leq 1,
\end{equation}
where $\gamma$ and $b_{\gamma}$ are constants. In particular, $\gamma=1$ corresponds to hard sphere molecules and $\gamma=0$ to Maxwell molecules. 

The collision operator $\mathcal{Q}$ has collision invariants $1$, $v$, and $|v|^2$;  that is,
\begin{equation} \label{consv}
\int_{\mathbb{R}^3}\mathcal{Q}(f)\,\rd{v}= \int_{\mathbb{R}^3}\mathcal{Q}(f) v\,\rd{v}=\int_{\mathbb{R}^3}\mathcal{Q}(f)|v|^2\,\rd{v}=0.
\end{equation}
In addition, $\mathcal{Q}$ satisfies Boltzmann's $H$-theorem;  that is,
\begin{equation} 
\int_{\mathbb{R}^3}\mathcal{Q}(f)\ln f\,\rd{v}\leq 0,
\end{equation}
with equality if and only if $f$ takes on the form of a Maxwellian:
\begin{equation} \label{Max}
M(v)=\frac{\rho}{(2\pi T)^{\frac{3}{2}}}e^{-\frac{|v-u|^2}{2T}},
\end{equation}
where the density $\rho$, bulk velocity $u$, and temperature $T$ are given by
\begin{align}
\rho=\int_{\mathbb{R}^3}f\,\rd{v}, \quad  u=\frac{1}{\rho}\int_{\mathbb{R}^3}fv\,\rd{v}, \quad T=\frac{1}{3\rho}\int_{\mathbb{R}^3}f|v-u|^2\,\rd{v}.
\end{align}
This implies that in the homogeneous case, the entropy $\mathcal{S}(f) = -\int_{\mathbb{R}^3} f \ln f\,\rd{v}$ is always non-decreasing and reaches its maximum at the equilibrium defined by the Maxwellian in \eqref{Max}.

Proposed by Ludwig Boltzmann in 1872, the Boltzmann equation (\ref{CBE}) is one of the fundamental equations of kinetic theory. Yet its numerical approximation still presents a huge computational challenge, even on today's supercomputers. This is mainly due to the high-dimensional, nonlinear, nonlocal structure of the collision integral in (\ref{CO}). Two approaches have been primarily employed for solving the Boltzmann equation numerically: one stochastic and one deterministic. Direct simulation Monte Carlo (DSMC) methods \cite{Bird, Nanbu80, Caflisch98} have been historically popular because they avoid the curse of dimensionality for this problem, however they can suffer from slow convergence for certain types of problems such as transient and low-speed flows and give noisy results due to their stochastic nature. The other approach is to use deterministic solvers, which have undergone considerable development over the past twenty years. These methods include discrete velocity models (DVM) \cite{RS94, BPS95, Buet96, MPR13} and Fourier spectral methods \cite{PP96, BR99, PR00, FR03, GT09, GT10}. DVMs are quadrature-based methods with grid points that are carefully chosen in order to preserve the conserved quantities of the collision operator. Spectral methods, on the other hand, compute the collision operator by exploiting its structure in Fourier space. Compared with DVM, they can provide significantly more accurate results with less numerical complexity; the conservation properties are not strictly maintained but are preserved up to spectral accuracy. Compared with DSMC, they produce smooth, noise free solutions and can simulate regimes that particle methods find difficult. 

Despite of the aforementioned advantages, spectral methods are invariably hindered in most real-world applications since they require  $O(N^6)$ operations per evaluation of the collision operator, with $N$ being the number of discretization points in each velocity dimension, as well as $O(N^6)$ bytes of memory to store precomputed weight functions.  With this type of scaling, the evaluation of the collision operator quickly becomes the bottleneck when solving large-scale problems \cite{PR00, GT09}. Fast spectral methods, based on the Carleman representation of the collision integral, have been proposed in \cite{BR99, MP06}.  These methods reduce the complexity of evaluating the collision operator to $O(MN^3\log N)$, where $M$ is the total number of discretization points on a sphere and $M\ll N^2$).  However, a decoupling assumption for the collision kernel is needed that restricts application of the method to the hard sphere case, i.e., $\gamma=1$ in (\ref{VHS}). In practice, however, $\gamma$ may take on any value in $[0,1]$; in addition, the collision kernel (\ref{CK}) may also have angular dependence. Therefore, the goal of this paper is to introduce a fast spectral method for the Boltzmann collision operator that can handle general collision kernels as well as mitigate the memory requirement in the direct spectral method. Specifically, the numerical complexity of our new method is $O(MN^4\log N)$; no precomputation is required for the VHS model, and only $O(MN^4)$ memory is needed to store precomputed functions for more general collision kernels. The proposed method can serve as a ``black-box" solver in the velocity domain to be used in conjunction with existing time and spatial discretization methods to treat more practical problems with complex geometries, multiple temporal/spatial scales, etc. Since our goal here is to present a simple strategy to accelerate the direct spectral method without sacrificing spectral accuracy, we will mainly focus on the approximation of the collision operator in the numerical examples and consider only the spatially homogeneous version of the Boltzmann equation \eqref{CBE}.

The rest of this paper is organized as follows. In the next section, we review the basic formulation of the direct spectral method and discuss its numerical challenges. The fast method is then described in Section 3. Numerical examples are presented in Section 4 to demonstrate the efficiency and accuracy of the proposed method. The paper is concluded in Section 5.

\section{The direct spectral method}

While multiple spectral formulations exist, we have elected in this paper to adopt the Fourier-Galerkin approach \cite{PR00} to illustrate the idea. The strategy introduced below can be easily applied to other spectral formulations such as the one based on the Fourier transform \cite{GT09}. 

The starting point for the spectral method is a change in the integration variable $v_*$ in (\ref{CO}) to $g=v-v_*$.  It is then assumed that $f$ has a compact support in $v$: $\text{Supp}_v(f)\approx \mathcal{B}_S$, where $\mathcal{B}_S$ is a ball centered at the origin with radius $S$. Of course, many distribution functions, including the Maxwellian, will not have compact support. Thus in practice, the support is chosen as some multiple (typically 6 to 8) of the thermal speed $v_{\rm{th}} = \sqrt{T}$.  It then suffices to truncate the infinite integral in $g$ to a larger ball $\mathcal{B}_R$ with radius $R=2S$:
\begin{equation}  \label{CO1}
\mathcal{Q}(f)(v)\approx \mathcal{Q}_R(f)(v) 
= \int_{\mathcal{B}_R}\int_{S^2}B\left(r,\omega\cdot \hat{g}\right)
\left[f(v')f(v_*')-f(v)f(v-g)\right]\,\rd{\omega}\,\rd{g},
\end{equation}
where
\begin{equation}
v'=v-\frac{g}{2}+\frac{r}{2}\omega  \quad \text{and} \quad  v_*'=v-\frac{g}{2}-\frac{r}{2}\omega,
\end{equation}
with $r =|g|$ and $\hat{g} = g/|g|$ being the magnitude and direction of $g$, respectively. Next, one restricts $f$ to the computational domain $\mathcal{D}_L=[-L,L]^3$ and extends it periodically to the whole space. For anti-aliasing purposes, we let $L\geq (3+\sqrt{2})S/2$.%
\footnote{See \cite{PR00} for more details justifying this choice of $L$.}  
Then $f$ is approximated by a truncated Fourier series:%
\footnote{Note that $k=(k_1,k_2,k_3)\in \mathbb{Z}^3$ is a multi-dimensional index so that, for example, the summation in \eqref{Fseries} is understood to be over the lattice $\{k\in \mathbb{Z}^3: -\frac{N}{2} \leq k_1,k_2,k_3 \leq \frac{N}{2}-1 \}$.}
\begin{equation} \label{Fseries}
f(v)\approx f_N(v)= \sum_{k={-\frac{N}{2}}}^{\frac{N}{2}-1}\hat{f}_k e^{i \frac{\pi}{L}k\cdot v}, \qquad \hat{f}_k=\frac{1}{(2L)^3}\int_{\mathcal{D}_L}f(v)e^{-i \frac{\pi}{L}k\cdot v}\,\rd{v}.
\end{equation}
By substituting (\ref{Fseries}) into (\ref{CO1}) and performing a Galerkin projection, one can express the $k$-th mode of the Fourier expansion of the collision operator as a (discrete) weighted convolution:
\begin{equation} \label{WC}
\hat{\mathcal{Q}}_k:=\frac{1}{(2L)^3} \int_{\mathcal{D}_L}\mathcal{Q}_R(f)(v) e^{-i \frac{\pi}{L}k\cdot v}\,\rd{v}
=\sum_{\substack{l,m=-\frac{N}{2}\\l+m=k}}^{\frac{N}{2}-1}\mathcal{G}(l,m)\hat{f}_l\hat{f}_m, \quad k=-\frac{N}{2},\dots,\frac{N}{2}-1,
\end{equation}
where $\mathcal{G}(l,m)=G(l,m)-G(m,m)$ and
\begin{equation} \label{weight1}
G(l,m)=\int_{\mathcal{B}_R}\int_{S^2}B\left(r,\omega\cdot \hat{g}\right)e^{-i\frac{\pi}{L}\frac{(l+m)}{2}\cdot g+i\frac{\pi}{L}r\frac{(l-m)}{2}\cdot \omega}\,\rd{\omega}\,\rd{g}.
\end{equation}
For the VHS model (\ref{VHS}), the formula for $G$ reduces to
\begin{align} 
G(l,m)= 16\pi^2b_{\gamma}\int_0^Rr^{\gamma+2}\,\text{Sinc}\left(\frac{\pi}{L}r\frac{|l+m|}{2}\right)\,\text{Sinc}\left(\frac{\pi}{L}r\frac{|l-m|}{2}\right)\,\rd{r},
\end{align} 
where $\text{Sinc}(x)=\sin x/x$.

To summarize, a single evaluation of the collision operator $\mathcal{Q}$ in the direct spectral method consists of the following steps:
\begin{enumerate}
\item[0.] precompute the weight $\mathcal{G}(l,m)$ --- storage requirement $O(N^6)$;
\item[1.] compute $\hat{f}_k$ using the fast Fourier transform (FFT) --- cost $O(N^3\log N)$;
\item[2.] compute the weighted convolution (\ref{WC}) --- cost $O(N^6)$;
\item[3.] take the inverse Fourier transform of $\hat{\mathcal{Q}}_k$ to obtain $\mathcal{Q}$ --- cost $O(N^3\log N)$.
\end{enumerate}
Step 2 is by far the most expensive step.  Indeed, due to the presence of the weights $\mathcal{G}({l,m})$ in the convolution, typical fast methods for convolutions do not apply.  Thus the constrained double summation in (\ref{WC}) has to be evaluated directly for every index $k$, resulting in $O(N^6)$ complexity.  Step 0 can be completed in advance, but it requires a huge amount of memory to store the precomputed weights.  This can be a challenge for large-scale problems, even on the largest supercomputers. For example, when $N=40$, it takes just over 32 gigabytes of data to store the weights---an amount that exceeds the memory capacity on a typical compute node of Oak Ridge National Laboratory's Titan supercomputer. 

\section{The new fast spectral method}

In the new fast spectral method, we accelerate the summation in (\ref{WC}). Let $\hat{\mathcal{Q}}_k=\hat{\mathcal{Q}}^+_k-\hat{\mathcal{Q}}^-_k$, where
\begin{align}\label{Q+}
\hat{\mathcal{Q}}^+_k=\sum_{\substack{l,m=-\frac{N}{2}\\l+m=k}}^{\frac{N}{2}-1}G(l,m)\hat{f}_l\hat{f}_m
\quad \text{and} \quad
\hat{\mathcal{Q}}^-_k=\sum_{\substack{l,m=-\frac{N}{2}\\l+m=k}}^{\frac{N}{2}-1}G(m,m)\hat{f}_l\hat{f}_m. 
\end{align}
Because $G(m,m)$ depends only on $m$, the loss term $\hat{\mathcal{Q}}^-_k$ is actually a convolution of the functions $G(m,m)\hat{f}_m$ and $\hat{f}_l$.  It can therefore be computed efficiently by FFT in $O(N^3\log N)$ operations, since convolution in the Fourier domain becomes multiplication in the original domain. What makes the total cost $O(N^6)$ is the gain term $\hat{\mathcal{Q}}^+_k$. 

Our goal is to find an approximation for $\hat{\mathcal{Q}}^+_k$ that can be expressed as a convolution. To this end, we seek an approximation of $G(l,m)$ in the following decoupled form:
\begin{equation} \label{lowrank}
G(l,m)\approx\sum_{p=1}^{N_p} \alpha_p(l+m)\beta_p(l)\gamma_p(m),
\end{equation}
where $\alpha_p$, $\beta_p$, and $\gamma_p$ are functions of $(l+m)$, $l$, and $m$ respectively, and $N_p \ll N^3$. Substitution of (\ref{lowrank}) into (\ref{Q+}) gives
\begin{equation} \label{sum}
\hat{\mathcal{Q}}^+_k\approx\sum_{\substack{l,m=-\frac{N}{2}\\l+m=k}}^{\frac{N}{2}-1}\sum_{p=1}^{N_p} \alpha_p(l+m)\beta_p(l)\gamma_p(m)\hat{f}_l\hat{f}_m=\sum_{p=1}^{N_p}\alpha_p(k)\sum_{\substack{l,m=-\frac{N}{2}\\l+m=k}}^{\frac{N}{2}-1}\left(\beta_p(l)\hat{f}_l\right)\left(\gamma_p(m)\hat{f}_m\right).
\end{equation}
The inner summation in (\ref{sum}) is a pure convolution of the two functions $\beta_p(l)\hat{f}_l$ and $\gamma_p(m)\hat{f}_m$ that can be computed in $O(N^3\log N)$ via FFT.  This, together with the outer summation, results in a total cost of $O(N_pN^3\log N)$ for a single evaluation of $\mathcal{Q}^+$. 

To generate a suitable low-rank approximation of the form (\ref{lowrank}), we propose a simple solution in which, instead of precomputing all the weights $G(l,m)$ in (\ref{weight1}), we compute them partially ``on the fly" using a quadrature rule. Specifically, we rewrite $G(l,m)$ as 
\begin{equation} \label{weight2}
G(l,m)=\int_0^R\int_{S^2}F(l+m,r,\omega)\,
e^{i\frac{\pi}{L}r\frac{l}{2}\cdot \omega}\,
e^{-i\frac{\pi}{L}r\frac{m}{2}\cdot \omega}\,
\rd{\omega}\,\rd{r},
\end{equation}
where
\begin{equation} 
F(l+m,r,\omega)=r^2\int_{S^2}B\left(r,\omega\cdot \hat{g}\right)e^{-i\frac{\pi}{L}r\frac{(l+m)}{2}\cdot \hat{g}}\,\rd{\hat{g}}.
\end{equation}
For each fixed $r$ and $\omega$, the integrand in (\ref{weight2}) is a product of three functions:  one that depends on $(l+m)$, one that depends on $l$, and one that depends on $m$.  This is exactly the desired
form of (\ref{lowrank}). In order to maintain this structure, we carry out the integration in $r$ and $\omega$ using a fixed numerical quadrature.  This yields
\begin{equation} \label{form}
G(l,m)\approx \sum_{r, \phi_1, \phi_2}w_{r}w_{\phi_1}w_{\phi_2} \sin\phi_2 \,F(l+m,r,\omega)e^{i\frac{\pi}{L}r\frac{l}{2}\cdot \omega}e^{-i\frac{\pi}{L}r\frac{m}{2}\cdot \omega},
\end{equation}
where $\phi_1$ is the azimuthal angle, $\phi_2$ is the polar angle, and $w_{r}$, $w_{\phi_1}$ and $w_{\phi_2}$ represent the corresponding quadrature weights. Since the radial direction oscillates on the scale of $\frac{l-m}{2}$, the number of quadrature points in $r$ must be at least $O(N)$ in order to resolve this dimension. For the integration on the sphere, we anticipate that the total number $M$ of quadrature points needed is much less than $N^2$;  this is confirmed in our numerical results. Thus we are able to obtain an admissible decomposition (\ref{lowrank}) of $G(l,m)$ with $N_p =  O(MN) \ll N^3$. Substituting (\ref{form}) into (\ref{Q+}), we have
\begin{align} \label{form1}
\hat{\mathcal{Q}}^+_k&\approx\sum_{r, \phi_1, \phi_2}w_{r}w_{\phi_1}w_{\phi_2}  \sin\phi_2 \, F(k,r,\omega)\sum_{\substack{l,m=-\frac{N}{2}\\l+m=k}}^{\frac{N}{2}-1}\left[ e^{i\frac{\pi}{L} r\frac{l}{2}\cdot \omega}\hat{f}_l \right] \left[e^{-i\frac{\pi}{L}r\frac{m}{2}\cdot \omega}\hat{f}_m\right].
\end{align}

In summary, the proposed fast algorithm for a single evaluation of $\mathcal{Q}$ consists of the following steps:
\begin{enumerate}
\item[0.] precompute the weights $G(m,m)$ and $F(k,r,\omega)$ --- storage requirement $O(MN^4)$;
\item[1.] compute $\hat{f}_k$ using FFT --- cost $O(N^3\log N)$;
\item[2.] compute the loss term $\mathcal{Q}^-$ using FFT --- cost $O(N^3\log N)$;
\item[3.] compute the gain term $\mathcal{Q}^+$ based on (\ref{form1}) using FFT --- cost $O(MN^4\log N)$;
\item[4.] compute $\mathcal{Q}=\mathcal{Q}^+-\mathcal{Q}^-$ --- cost $O(N^3)$.
\end{enumerate}
Compared with the direct spectral method in the previous section, the new method saves both in memory storage (step 0) and computational time (step 3). For the $N=40$ case mentioned in Section 2, if we take $M=14$, the precomputed weights only require roughly $247$ megabytes.  This is less than one percent of the memory required for the direct method. For the VHS model, the function $F$ does not depend on $\omega$ and has an analytical form
\begin{equation} 
F(k,r)=4\pi \,b_{\gamma}\,r^{\gamma+2}\, \text{Sinc}\left(\frac{\pi}{L}r\frac{|k|}{2}\right).
\end{equation}
Thus no precomputation is needed in this case.

In our numerical implementation, we use the Gauss-Legendre quadrature in the radial direction $r$, while for the integration in $\omega$ we propose to use the Lebedev quadrature \cite{LL99} which is the near optimal quadrature on the sphere and requires fewer quadrature points than tensor product based Gauss quadratures for a large class of functions \cite{Beentjes15}. The Lebedev quadrature is designed to enforce the exact integration of spherical harmonics up to a given order and only a certain number of quadrature points are available. To gain a concrete idea of how many quadrature points $M$ is needed for our problem, note that in a typical numerical example where $N=32$, we only need $M=14$ to reach a relative $10^{-6}$ accuracy; a larger $M$ (say, $M=74$) may be needed when considering anisotropic distribution functions. Nevertheless, it is generally expected that $M\ll N^2$. A similar observation has been made for the fast spectral method in \cite{MP06, FMP06, Filbet12}. Although this method is based on a different representation of the collision integral (and restricted only to the hard sphere case), it also requires numerical discretization on a sphere.


\begin{remark}
The method proposed above can be followed by a post-processing subroutine after each evaluation of the collision operator to strictly enforce the collision invariants in (\ref{consv}) \cite{GT09} for either scalar or  system Boltzmann models, where it is shown that the solution of the scalar problem converges to the equilibrium Maxwellian state (\ref{Max}) \cite{AGT16}, or alternatively, adapted easily to preserve  exactly the equilibrium Maxwellian state as proposed in \cite{FPR15}. Since the goal in this paper is to present a simple strategy to accelerate the computation of the weighted convolution (\ref{WC}) in the direct spectral method without sacrificing spectral accuracy, we will mainly focus on the proposed method itself in the following numerical examples and leave the investigation of aforementioned extensions to future work.
\end{remark}


\section{Numerical examples}

In this section, we perform a series of numerical tests to validate the accuracy and efficiency of the proposed method. In the first test, we compare results of the new method to the Bobylev-Krook-Wu (BKW) solution \cite{Bobylev75_1, KW77}, which is constructed for Maxwell-type interactions (i.e., $\gamma=0$ in (\ref{VHS})) and is one of the few analytical solutions available for the Boltzmann equation.  In the second test, we again consider Maxwell molecules, but assume an initial condition that is anisotropic in $v$.  In this case, there is no analytical solution for the full distribution function, but exact formulas for higher order moments such as the momentum flow tensor
\begin{equation}
\label{P}
	P_{ij}=\int_{\mathbb{R}^3} f v_i v_j\,\rd{v}, \quad i,j=1,2,3,
\end{equation}
	 and the energy flow vector
\begin{equation}
\label{q}
q_{i}=\frac{1}{2}\int_{\mathbb{R}^3} f v_i v^2\,\rd{v}, \quad i=1,2,3
\end{equation}
can be derived.  We will also test our method for the hard sphere case by comparing it with the direct spectral method. Finally, we illustrate the generality of our method by considering a more realistic, angularly dependent collision kernel.

In the following, ``direct spectral" refers to the direct spectral method, and ``fast spectral" refers to the new method proposed in this paper. The implementation is in MATLAB and all numerical results are obtained on a laptop computer (MacBook Pro, 3.0GHz Dual-core Intel Core i7 with 8GB memory). Further acceleration can be achieved by careful implementation in C or Fortran.

\subsection{Maxwell molecules -- BKW solution}

When the collision kernel $B=1/(4\pi)$ is a constant, one can construct an exact solution to the spatially homogeneous Boltzmann equation
\begin{equation} \label{HCBE}
\frac{\partial f}{\partial t}=\mathcal{Q}(f), \quad t>0, \ v\in \mathbb{R}^3.
\end{equation}
This solution takes the form
\begin{equation} \label{ext}
f(t,v)=\frac{1}{2(2\pi K(t))^{3/2}}
\exp\left(-\frac{v^2}{2K(t)}\right)\left(\frac{5K(t)-3}{K(t)}+\frac{1-K(t)}{K^2(t)}v^2\right),
\end{equation}
where $K(t)=1-\exp(-t/6)$. The initial time $t_0$ has to be greater than $6\ln (5/2) \approx 5.498$ for $f$ to be positive. We take $t_0=5.5$.

Since $f$ given in (\ref{ext}) satisfies (\ref{HCBE}) exactly, the time derivative of $f$ gives
\begin{multline} \label{extQ}
\mathcal{Q}(f) \equiv \frac{\partial f}{\partial t}
= \left\{\left( -\frac{3}{2K} +\frac{v^2}{2K^2}  \right) f
+ \frac{1}{2(2\pi K)^{3/2}}\exp\left(-\frac{v^2}{2K}\right) 
\left( \frac{3}{K^2}+\frac{K-2}{K^3}v^2\right)\right\}K',
\end{multline}
where $K'(t)=\exp(-t/6)/6$. Using (\ref{extQ}), we can verify the accuracy of the proposed method without introducing additional time discretization error. We pick an arbitrary time $t=6.5$ and compare the numerical error and evaluation time of the direct and fast spectral methods. The results are shown in Tables \ref{table3} and \ref{table4}, from which we see that only 14 points are needed on the sphere for the fast method to obtain comparable accuracy to the direct method. Meanwhile, the speedup is about a factor of 300, even for the moderate value $N=32$. When $N=64$, the direct method requires too much storage for precomputed weights to fit within the available memory; it is therefore omitted. This restriction highlights the advantage of the proposed method in terms of memory.

\begin{table}[h!]
\centering
 \begin{tabular}{ c | c | c }
\hline 
 $N$ & direct spectral & fast spectral $M=14$ \\
 \hline
 $8$ & 6.91e-04 &7.33e-04   \\
  $16$&7.83e-05 &7.63e-05 \\
   $32$& 3.90e-08 &3.90e-08   \\
  $64$ & ---  & 3.81e-08  \\
 \hline
\end{tabular}
\caption{$\|\mathcal{Q}^{\text{num}}(f)-\mathcal{Q}^{\text{ext}}(f)\|_{L^{\infty}}$ evaluated at $t=6.5$. $N$ is the number of discretization points in each velocity dimension. In the fast spectral method, $N$-point Gauss quadrature is used in the radial direction and $M=14$ Lebedev rule is used for the sphere integration. We set $R=6$ (integration range), and $L=(3+\sqrt{2})R/4\approx6.62$ (computational domain).}
\label{table3}
\end{table}

\begin{table}[h!]
\centering
 \begin{tabular}{ c | c | c   }
\hline 
 $N$ & direct spectral & fast spectral $M=14$ \\
 \hline
 $8$ & 0.09s  & 0.14s \\
  $16$& 6.31s &0.26s \\
   $32$& 542.34s & 1.78s \\
  $64$ & --- &33.15s   \\
 \hline
\end{tabular}
\caption{Average running time for one time evaluation of the collision operator. Same parameters as in Table \ref{table3}.}
\label{table4}
\end{table}

\begin{figure}[htp]
\begin{center}
\includegraphics[width=3.5in,height=2.5in]{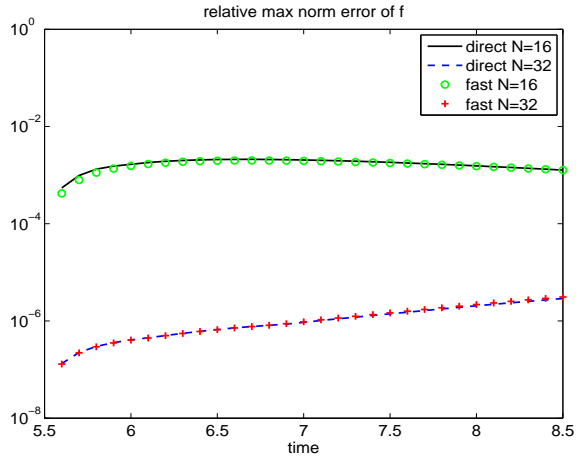}
\end{center}
\caption{Time evolution of $\|f^{\text{num}}-f^{\text{ext}}\|_{L^{\infty}}/\|f^{\text{ext}}\|_{L^{\infty}}$. RK4 with $\Delta t=0.1$ for time discretization. Other parameters are the same as in Table \ref{table3}.}
\label{figure2}
\end{figure}

We next couple the fast collision solver with time discretization to numerically solve the Boltzmann equation (\ref{HCBE}). A $4$th-order Runge-Kutta method is used to ensure that the temporal error does not pollute the spectral accuracy in velocity. The results are shown in Figure \ref{figure2}. The fast method basically produces very similar results to the direct method. Other norms behave similarly, but are omitted for brevity.
 
\subsection{Maxwell molecules -- moments}

Consider again the constant collision kernel $B=1/(4\pi)$. For the initial condition
\begin{equation} 
f(0,v)=\frac{1}{2(2\pi)^{3/2}}\left\{ \exp\left(-\frac{(v-u_1)^2}{2}\right)  + \exp\left(-\frac{(v-u_2)^2}{2}\right)  \right\},
\end{equation}
with $u_1=(-2,2,0)$ and $u_2=(2,0,0)$, the exact formulas for the non-zero components of $P$ and $q$ (cf. \eqref{P} and \eqref{q}) are given by
\begin{align}
&P_{11}=\frac{7}{3}\exp\left(-\frac{t}{2} \right)+\frac{8}{3}, \quad P_{22}=-\frac{2}{3}\exp\left(-\frac{t}{2} \right)+\frac{11}{3}, \nonumber\\
&P_{33}=-\frac{5}{3}\exp\left(-\frac{t}{2} \right)+\frac{8}{3}, \quad P_{12}=-2\exp\left(-\frac{t}{2} \right),
\end{align}
and
\begin{equation}
q_1=-2\exp\left(-\frac{t}{2} \right), \quad q_2=-\frac{2}{3}\exp\left(-\frac{t}{2} \right)+\frac{43}{6}.
\end{equation}

\begin{figure}[htp]
\begin{center}
\includegraphics[width=2.86in]{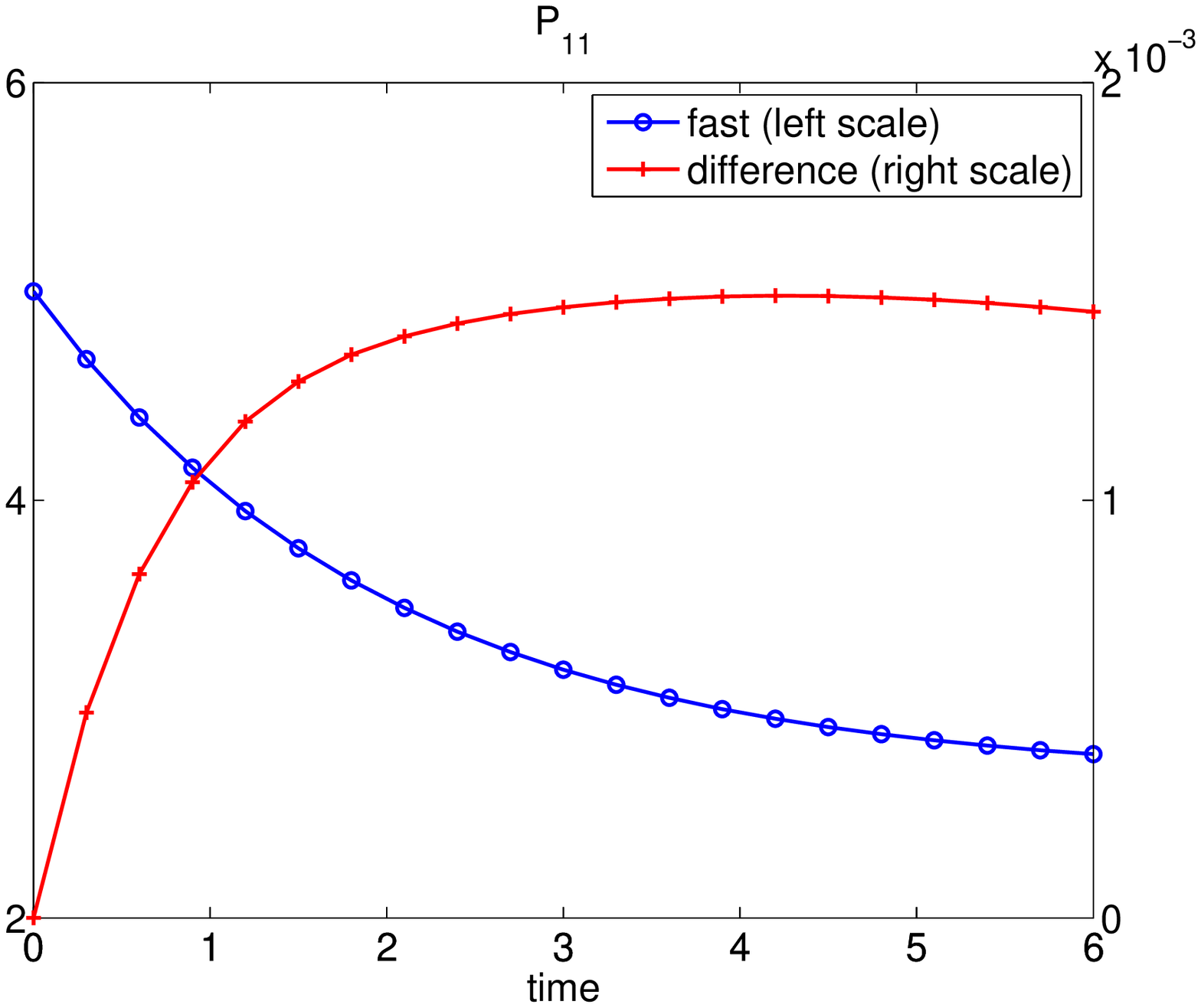}
\includegraphics[width=2.86in]{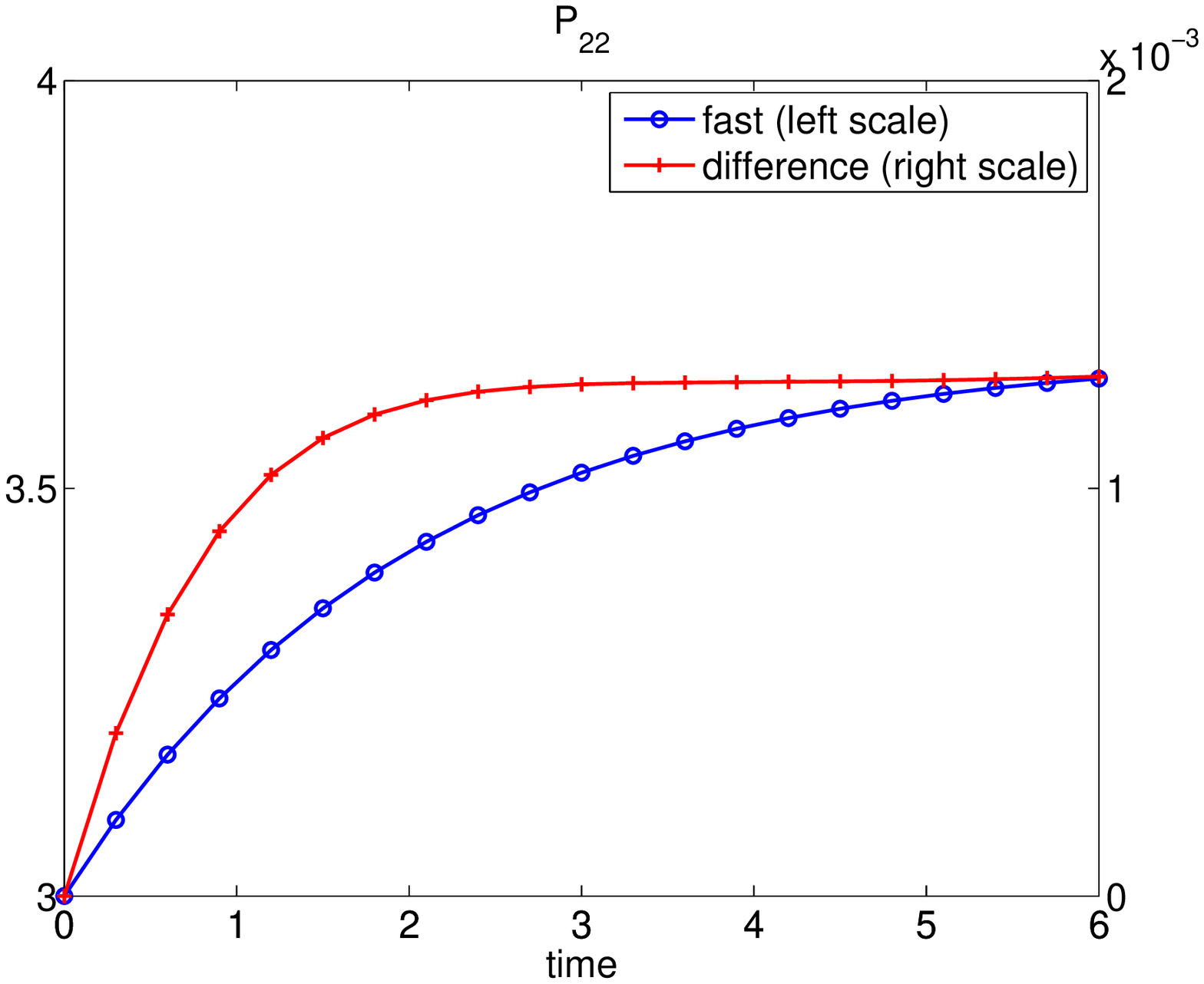}
\includegraphics[width=2.86in]{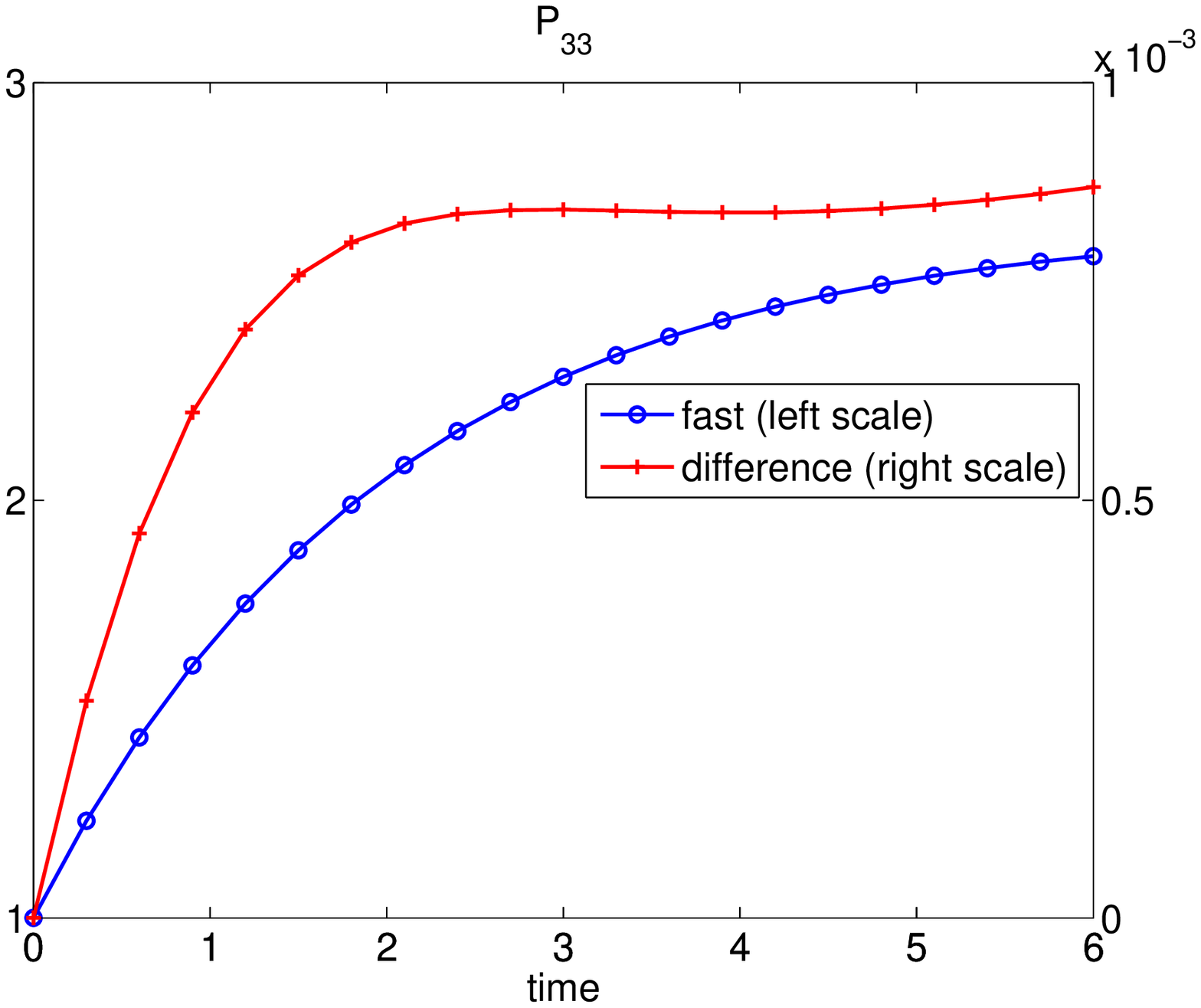}
\includegraphics[width=2.86in]{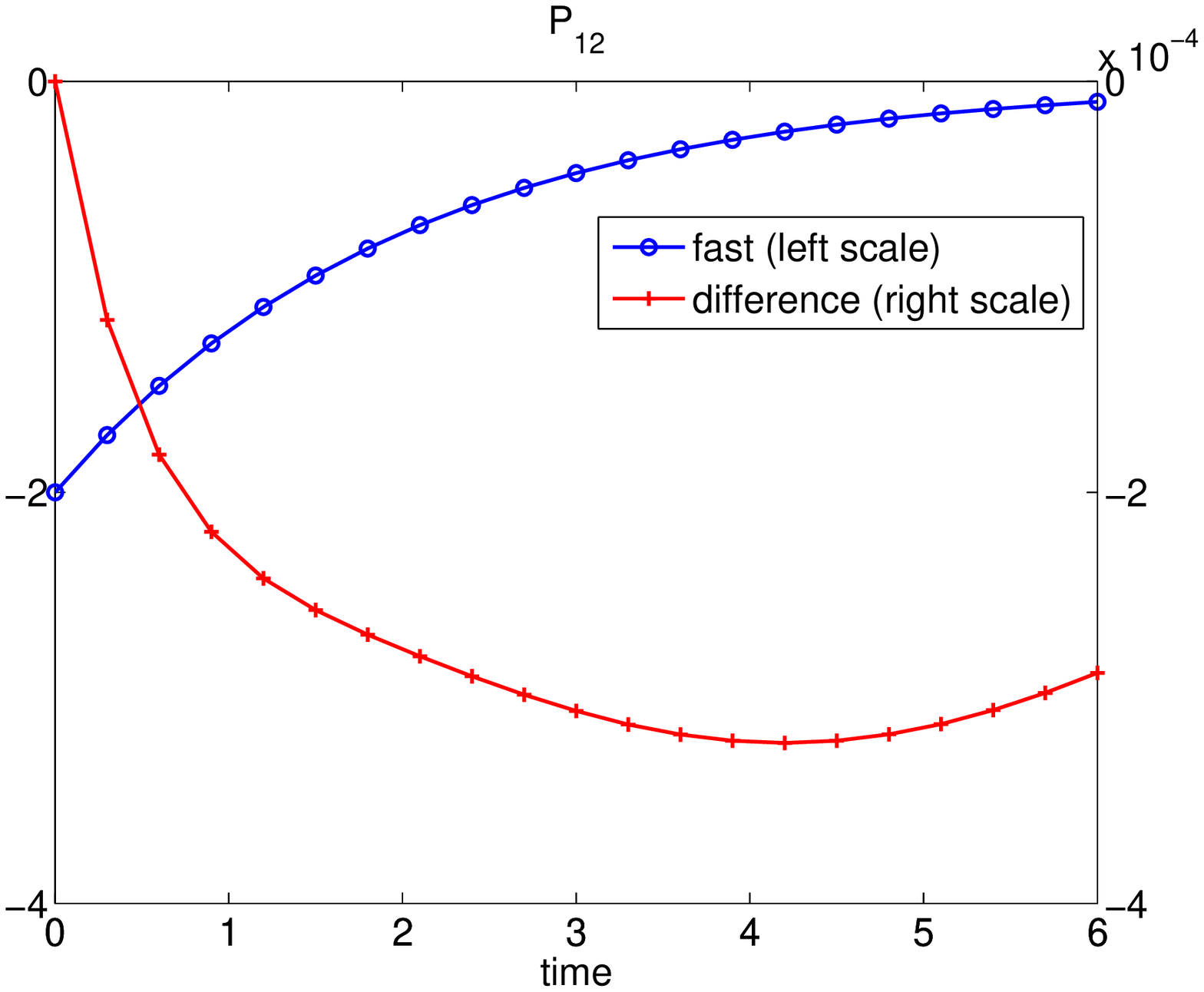}
\includegraphics[width=2.86in]{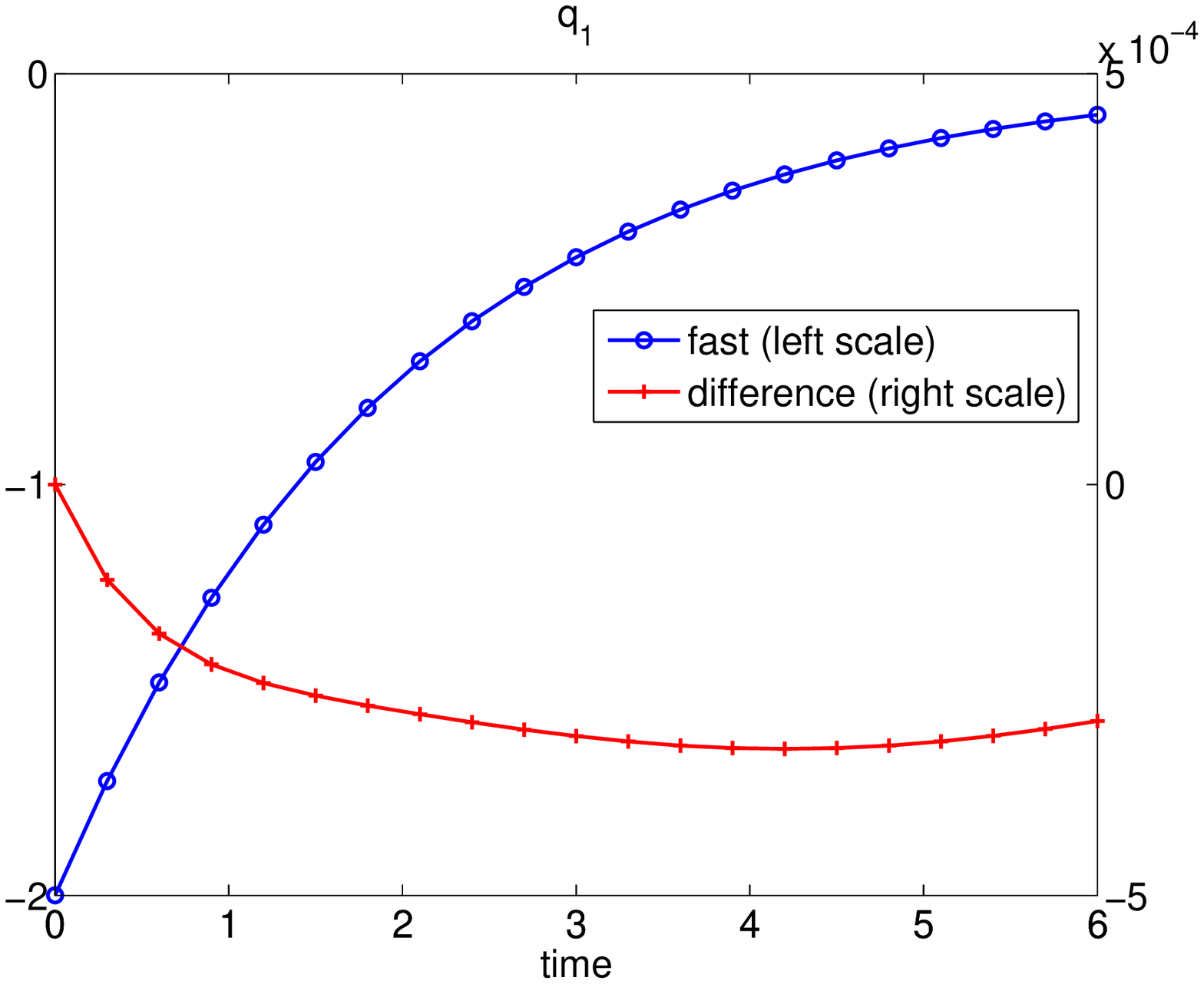}
\includegraphics[width=2.86in]{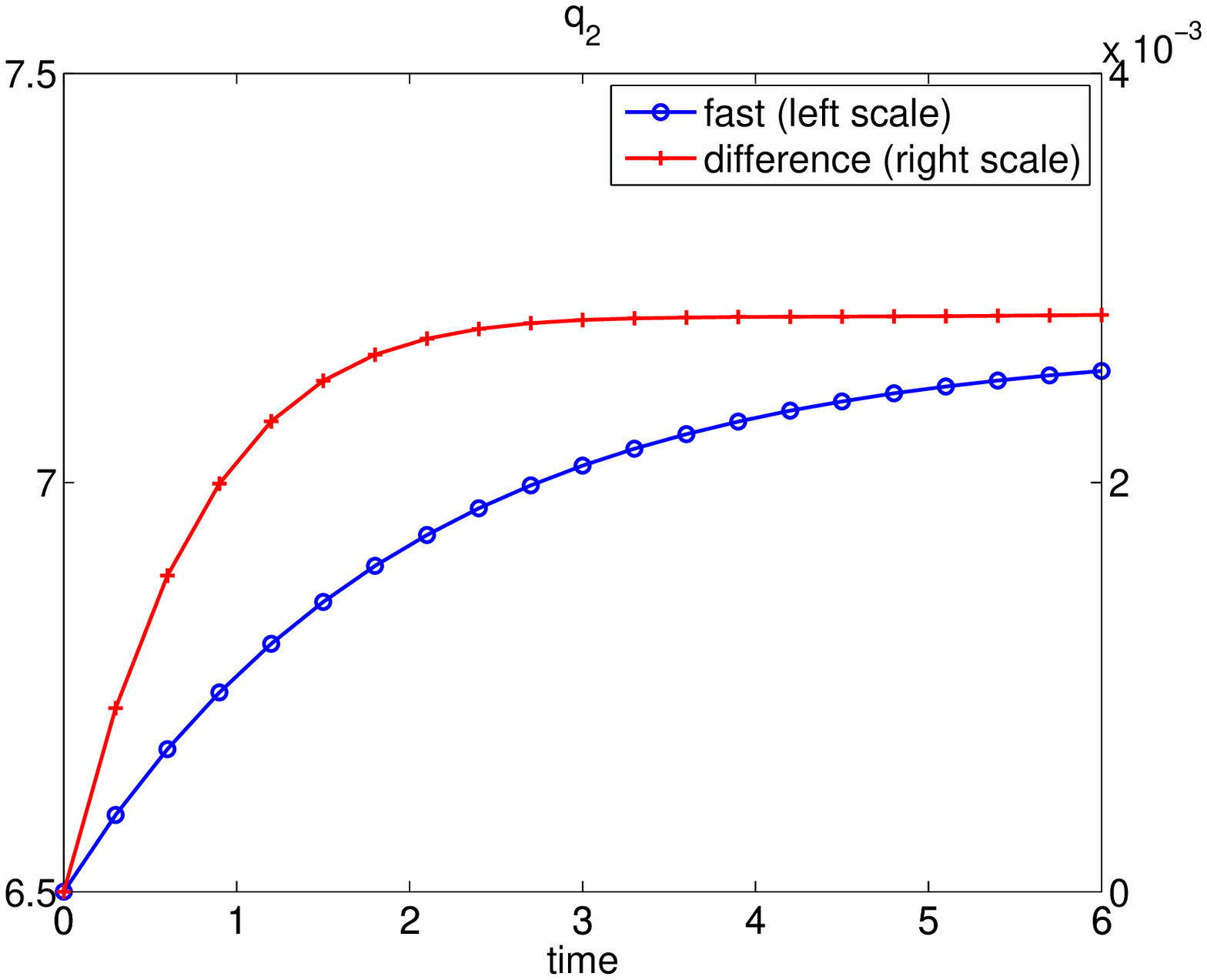}
\end{center}
\caption{Maxwell molecules. Time evolution for higher order moments. In each figure, the left scale shows the result by the fast spectral method, and the right scale shows its difference from the exact solution. RK4 with $\Delta t=0.3$ for time discretization. $N=32$ in each velocity dimension. In the fast method, $N=32$ in radial direction and $M=74$ for sphere integration. $R=10$, $L=(3+\sqrt{2})R/4\approx11.04$.}
\label{figure3}
\end{figure}

\begin{figure}[htp]
\begin{center}
\includegraphics[width=1.89in]{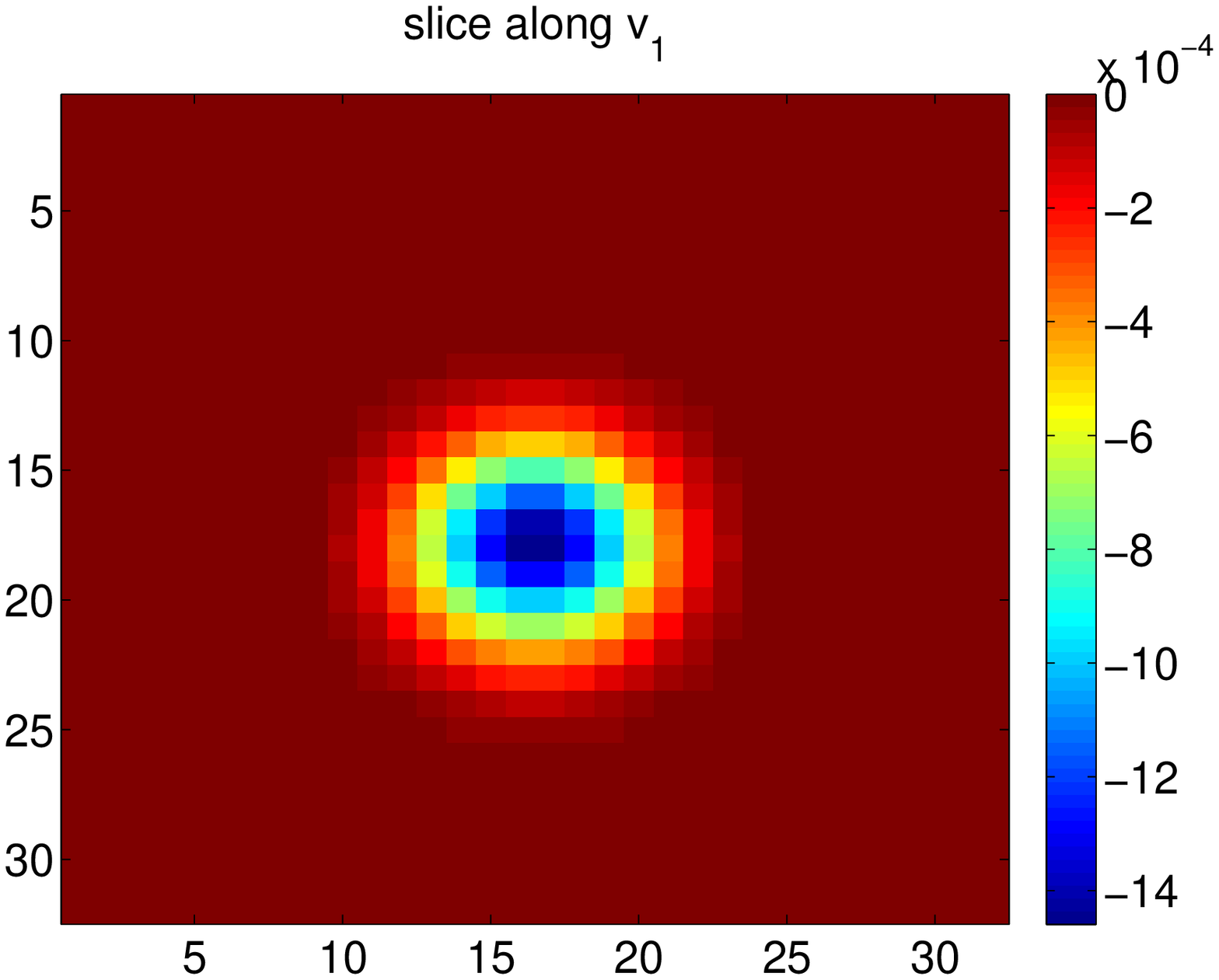}
\includegraphics[width=1.89in]{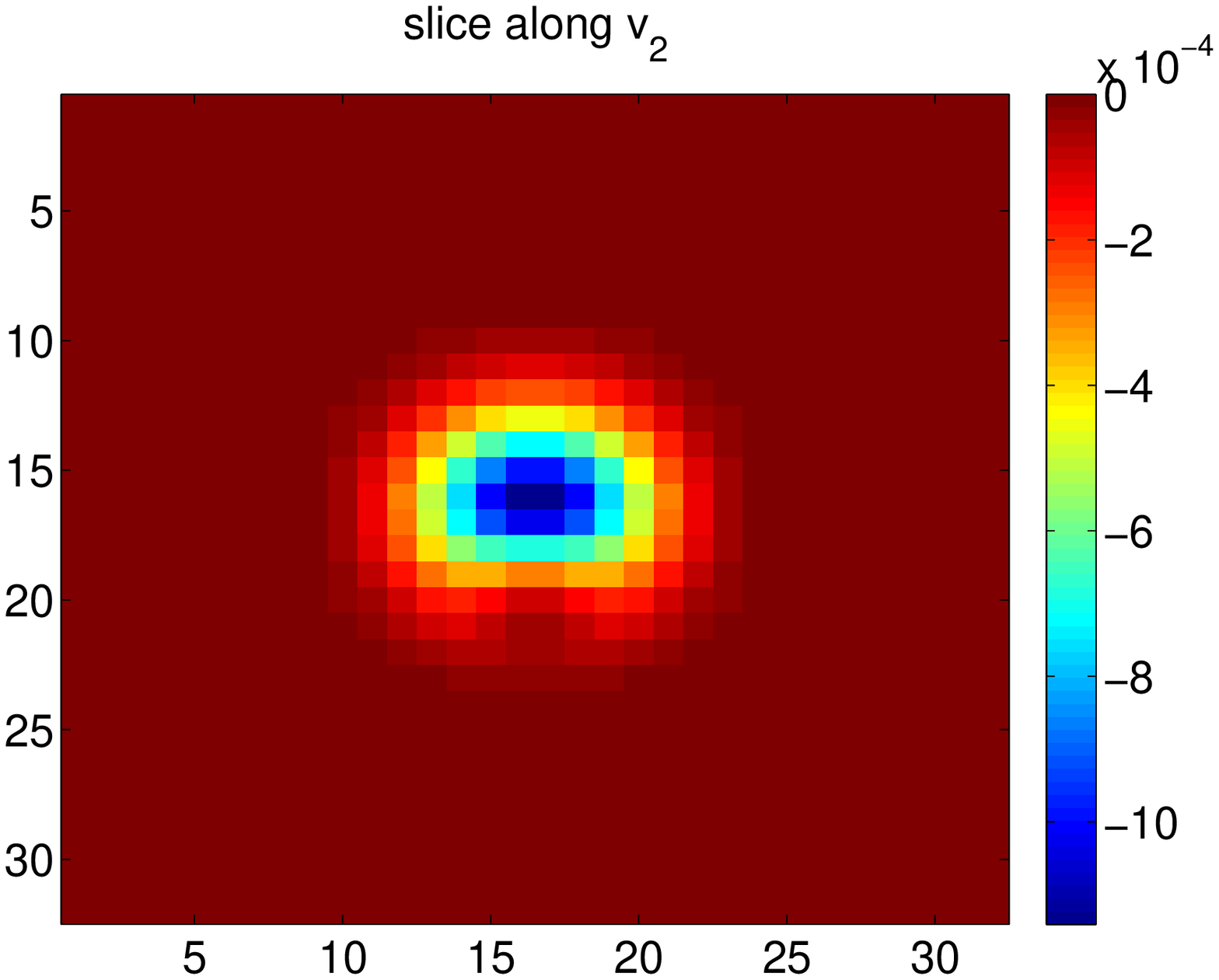}
\includegraphics[width=1.89in]{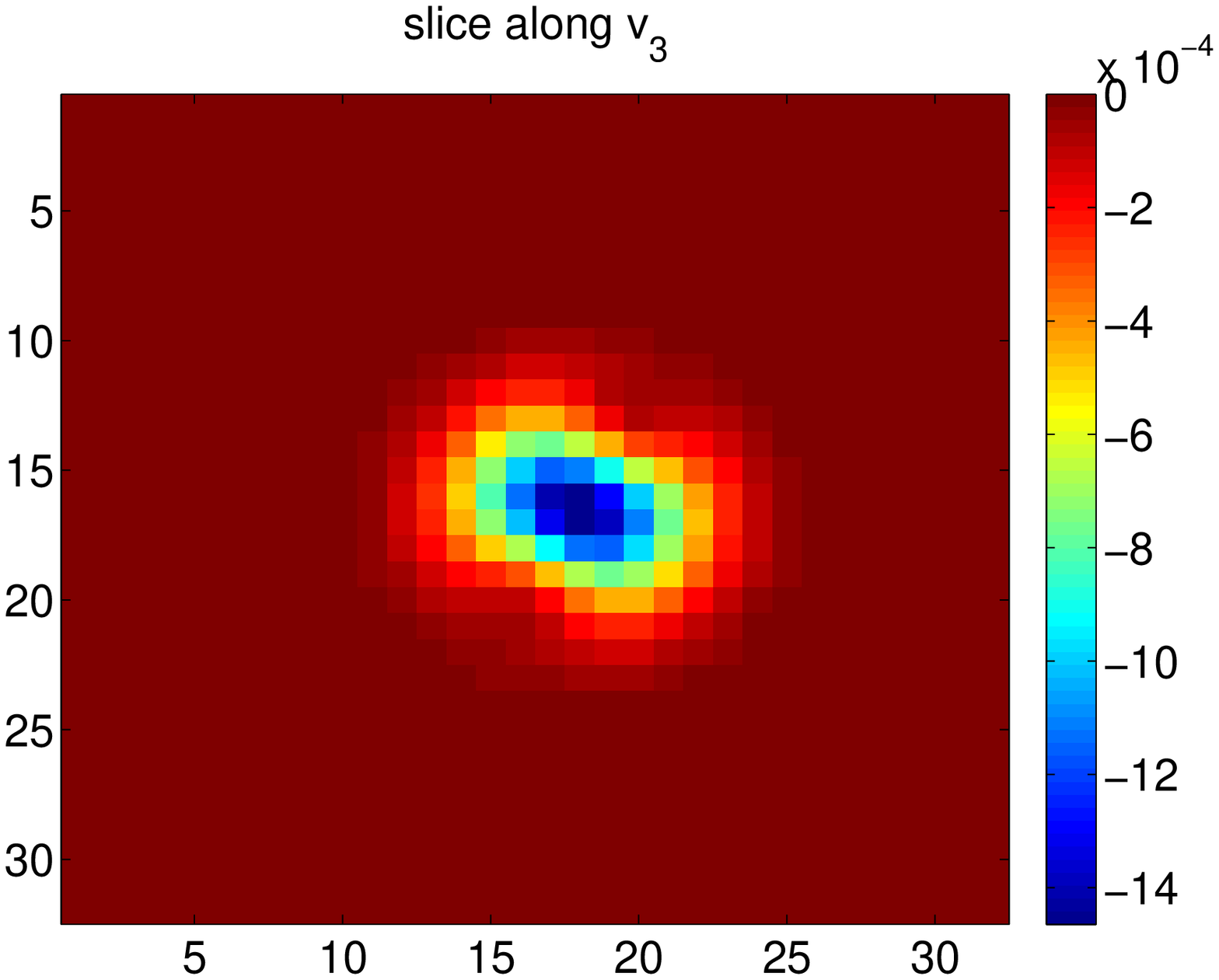}
\includegraphics[width=1.89in]{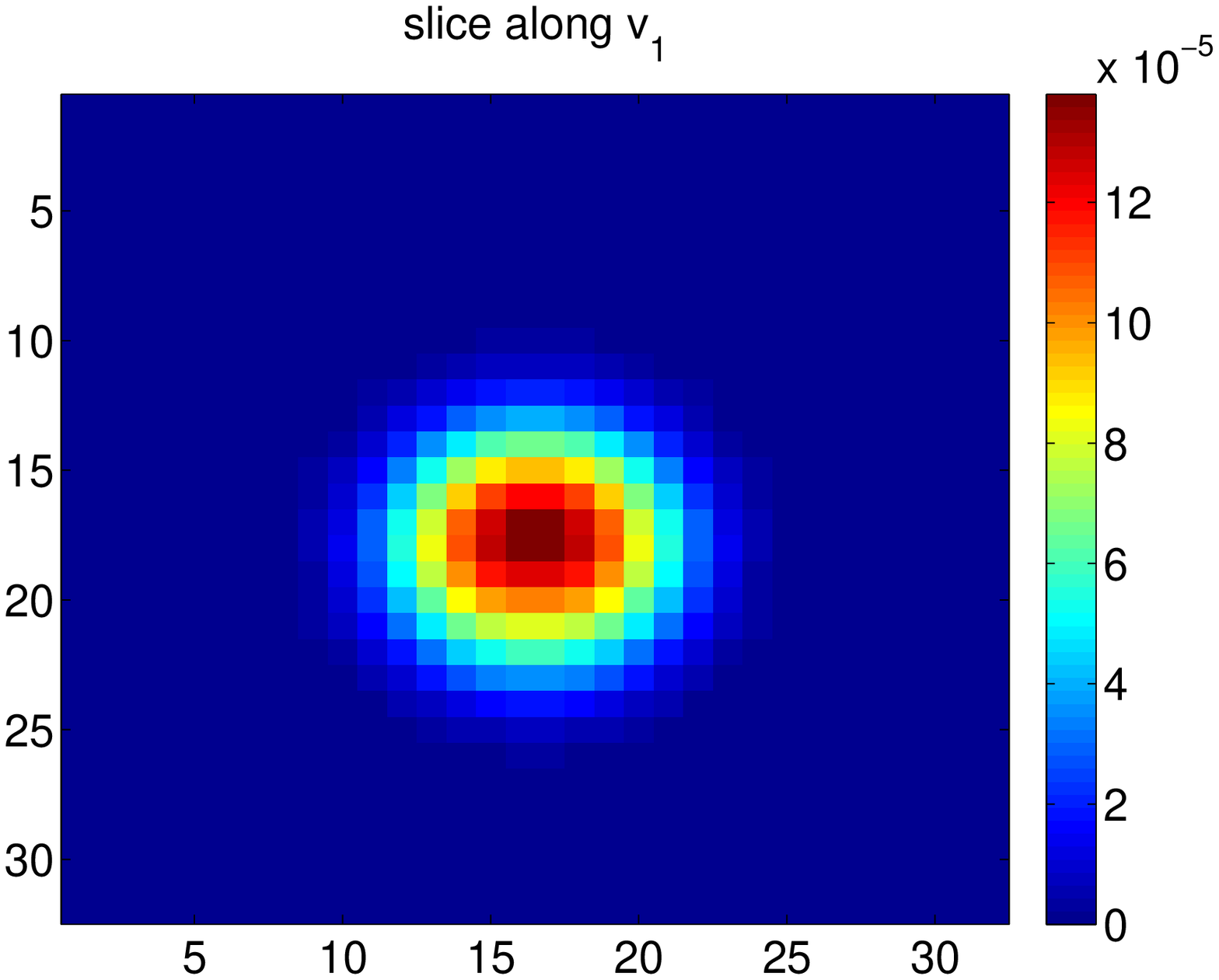}
\includegraphics[width=1.89in]{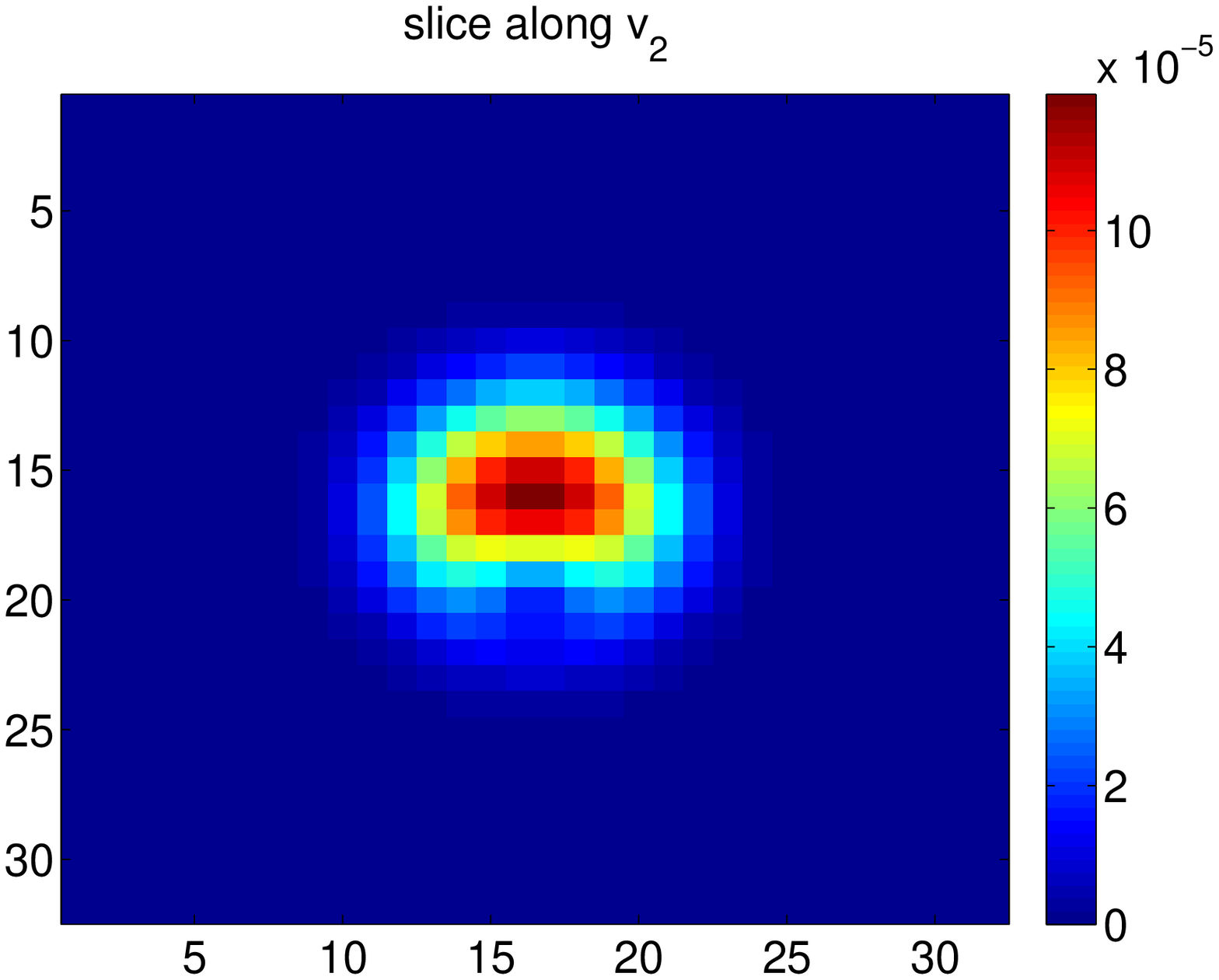}
\includegraphics[width=1.89in]{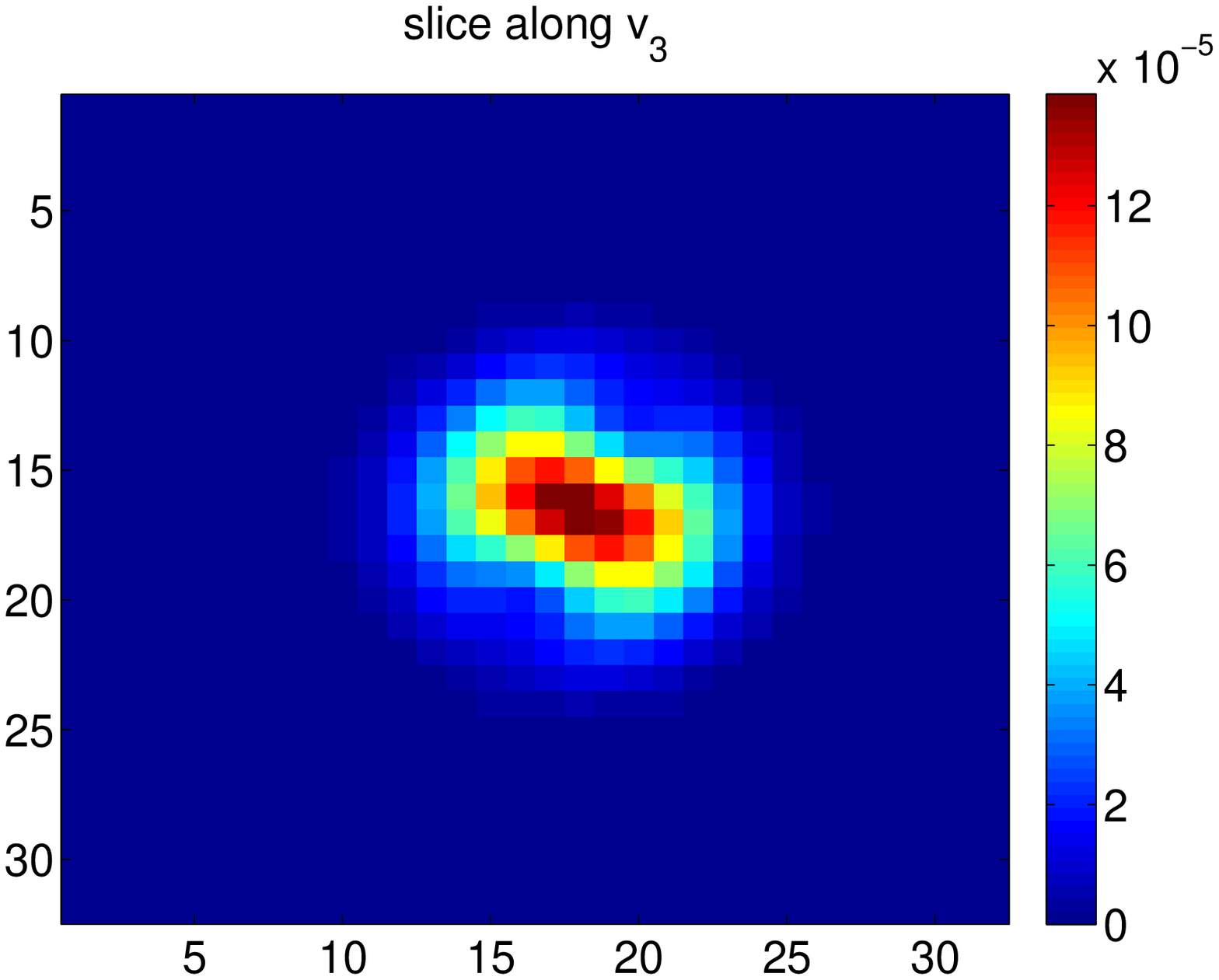}
\includegraphics[width=1.89in]{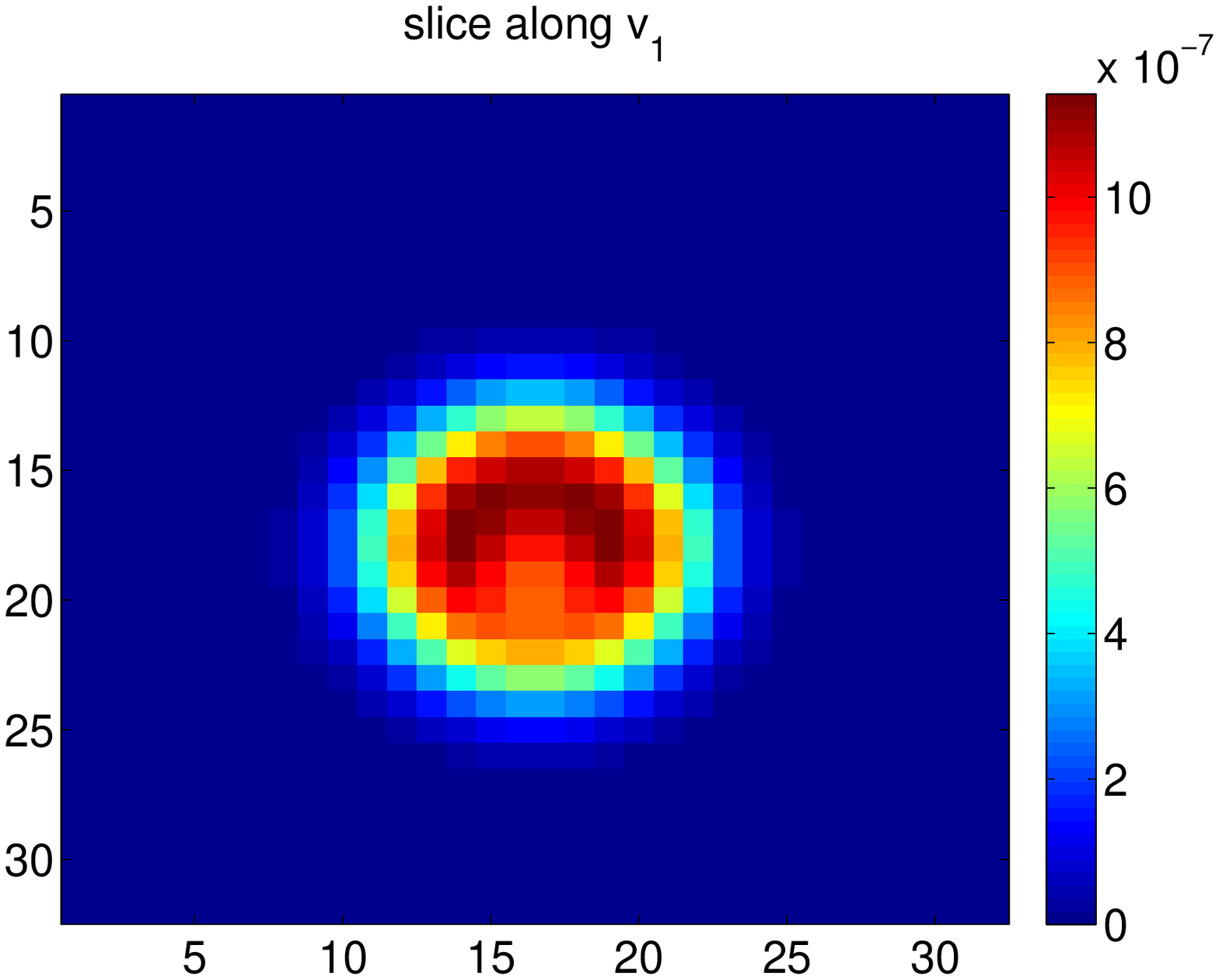}
\includegraphics[width=1.89in]{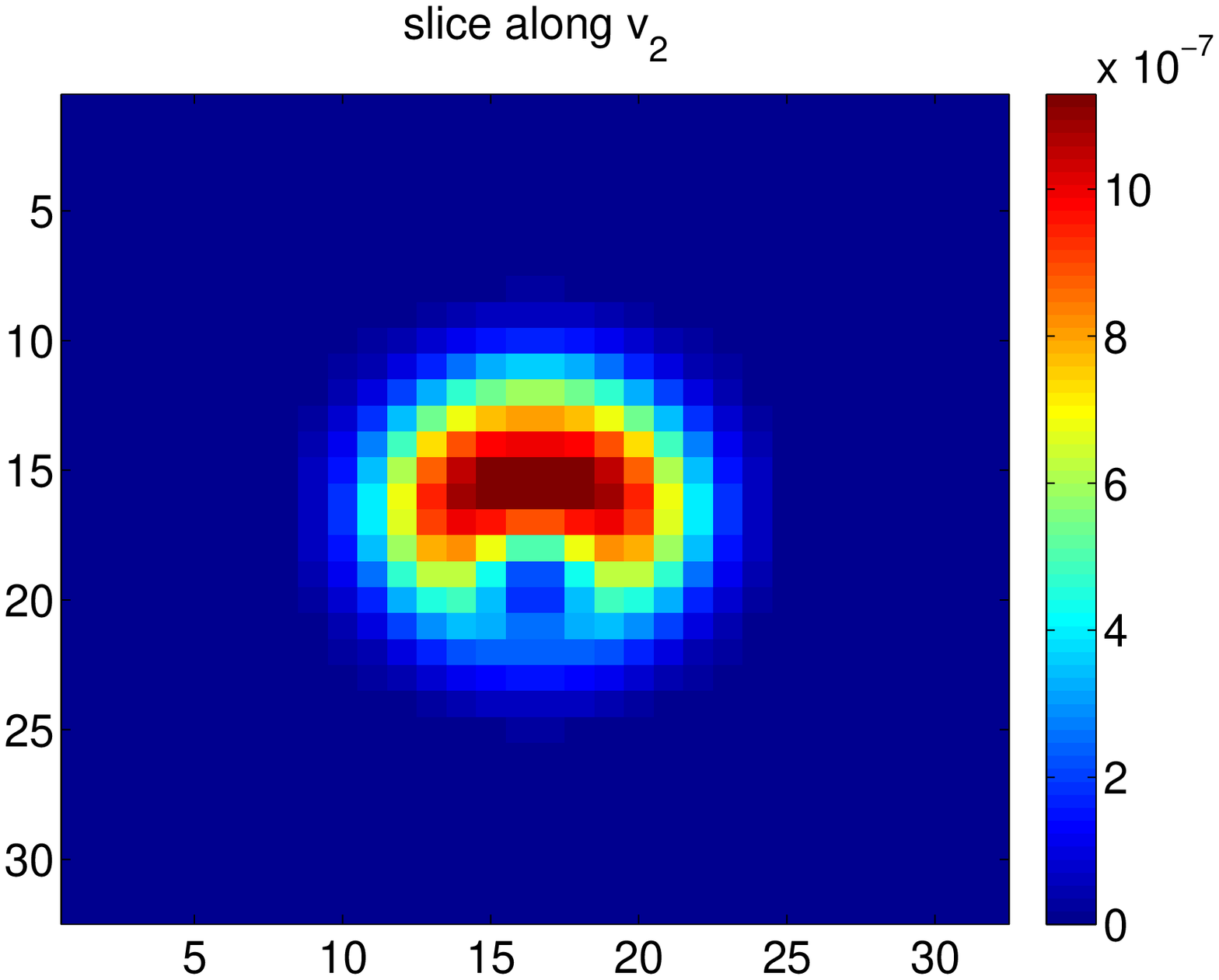}
\includegraphics[width=1.89in]{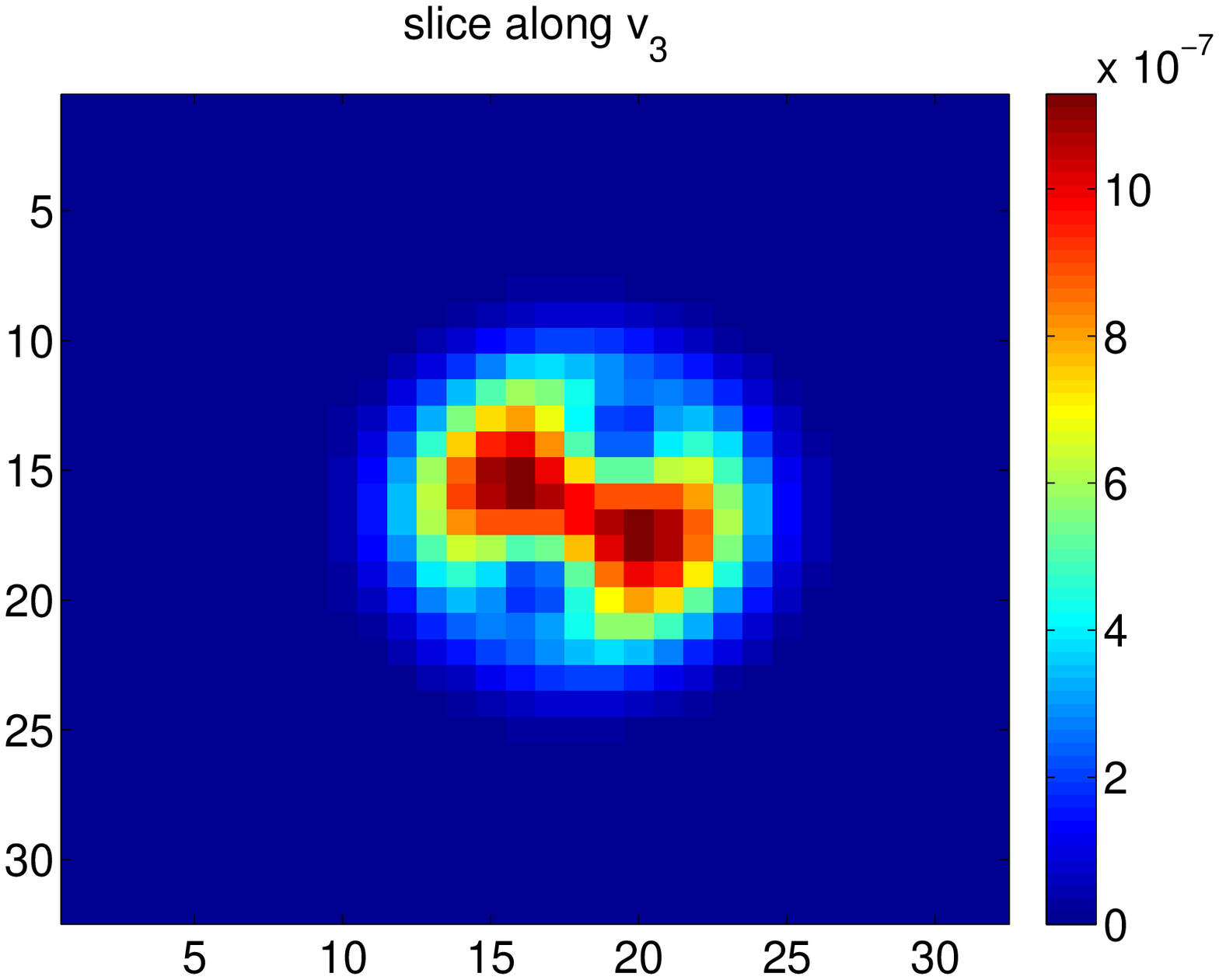}
\end{center}
\caption{$f^{\text{fast}}-f^{\text{direct}}$ at time $t=6$. From left to right: slices along the direction of $v_1$, $v_2$ and $v_3$. From top to bottom: $f^{\text{fast}}$ obtained with $M=14$, $38$ and $74$. Other parameters are the same as in Figure \ref{figure3}.}
\label{figure3_1}
\end{figure}

In Figure \ref{figure3}, we compare the results of the fast method with the formulas above. Because the solution is anisotropic in $v$, we need a larger value of $M$ than before to obtain reasonable accuracy. We find that $M=74$, which is still much less than $N^2=1024$, is sufficient to obtain roughly three digits of accuracy for moments. In Figure \ref{figure3_1}, we plot differences in the solution of the distribution function computed with the fast method and the direct method for different values of $M$.  (Since the exact distribution function is not known, we use the solution of the direct method as a reference.) We observe that the error decreases quickly as $M$ increases.

\subsection{Hard sphere molecules -- moments}

\begin{figure}[htp]
\begin{center}
\includegraphics[width=2.86in]{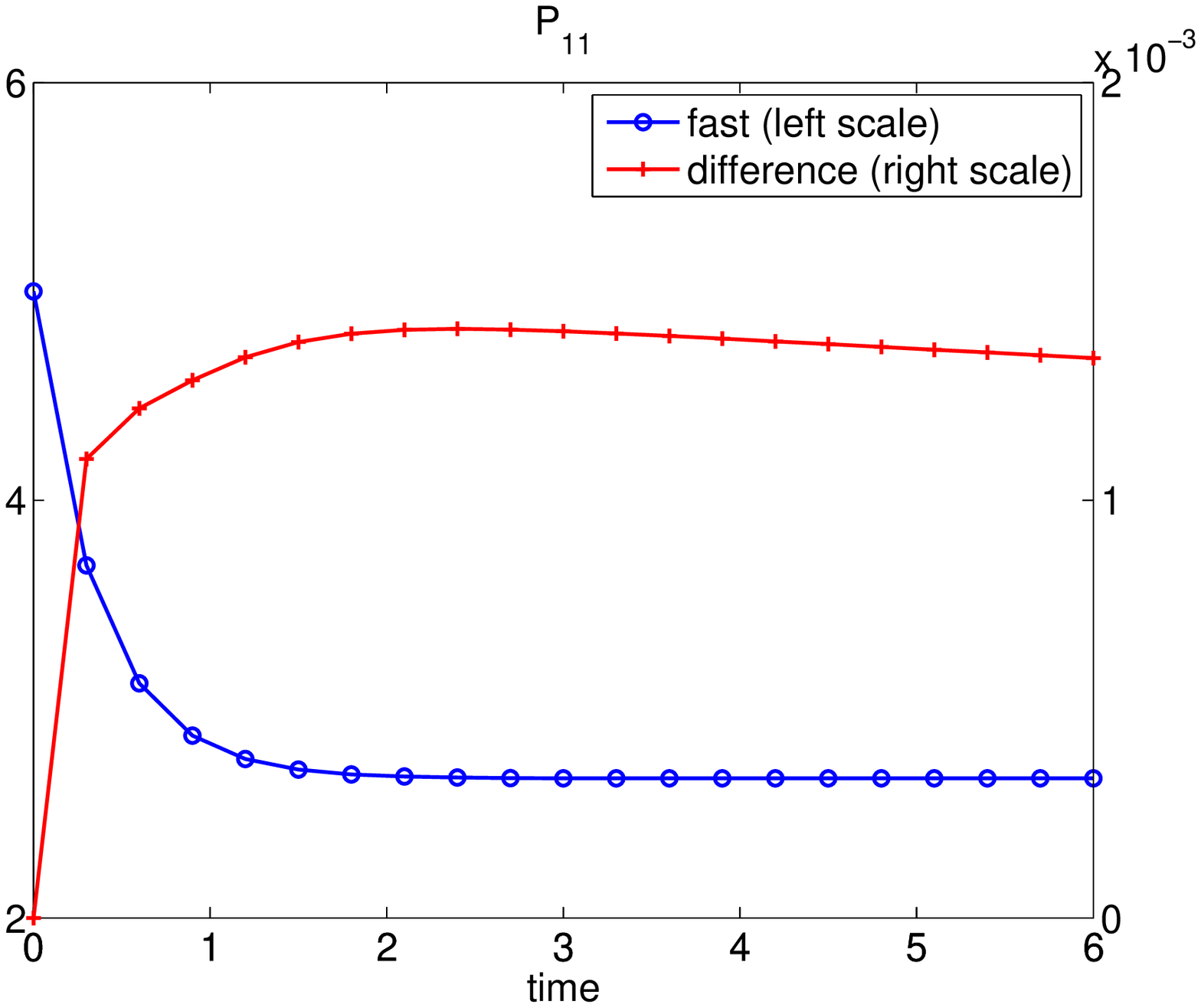}
\includegraphics[width=2.86in]{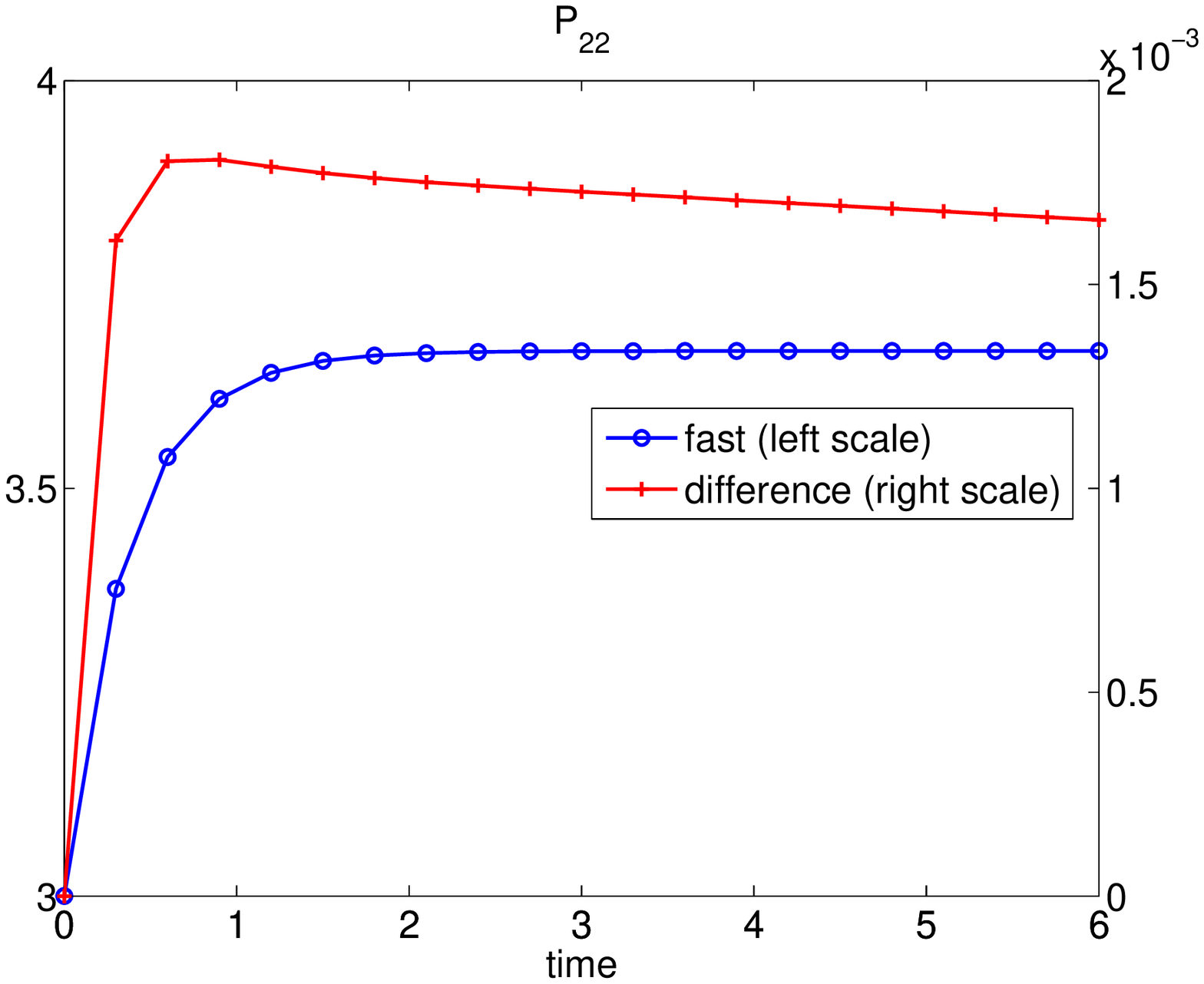}
\includegraphics[width=2.86in]{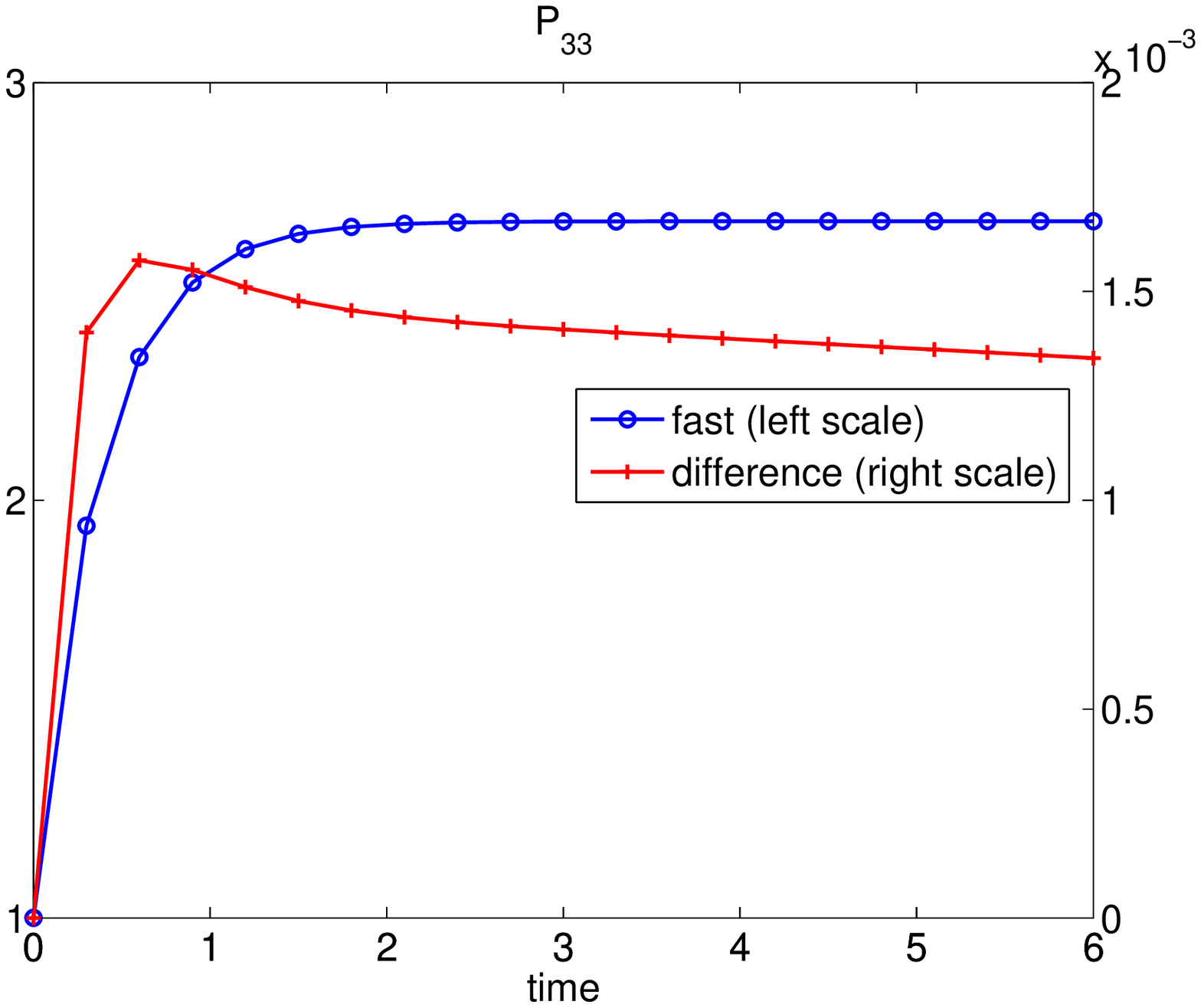}
\includegraphics[width=2.86in]{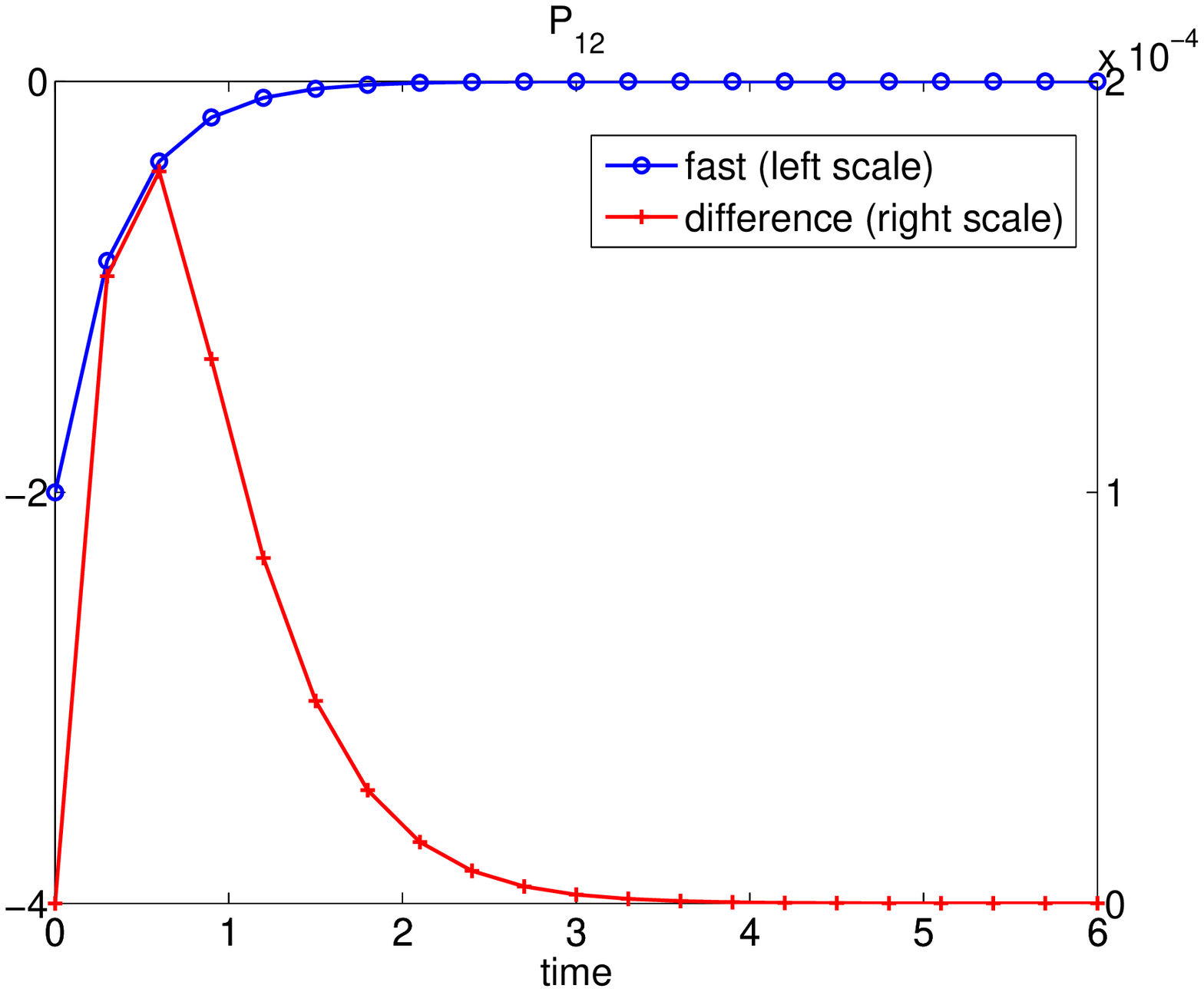}
\includegraphics[width=2.86in]{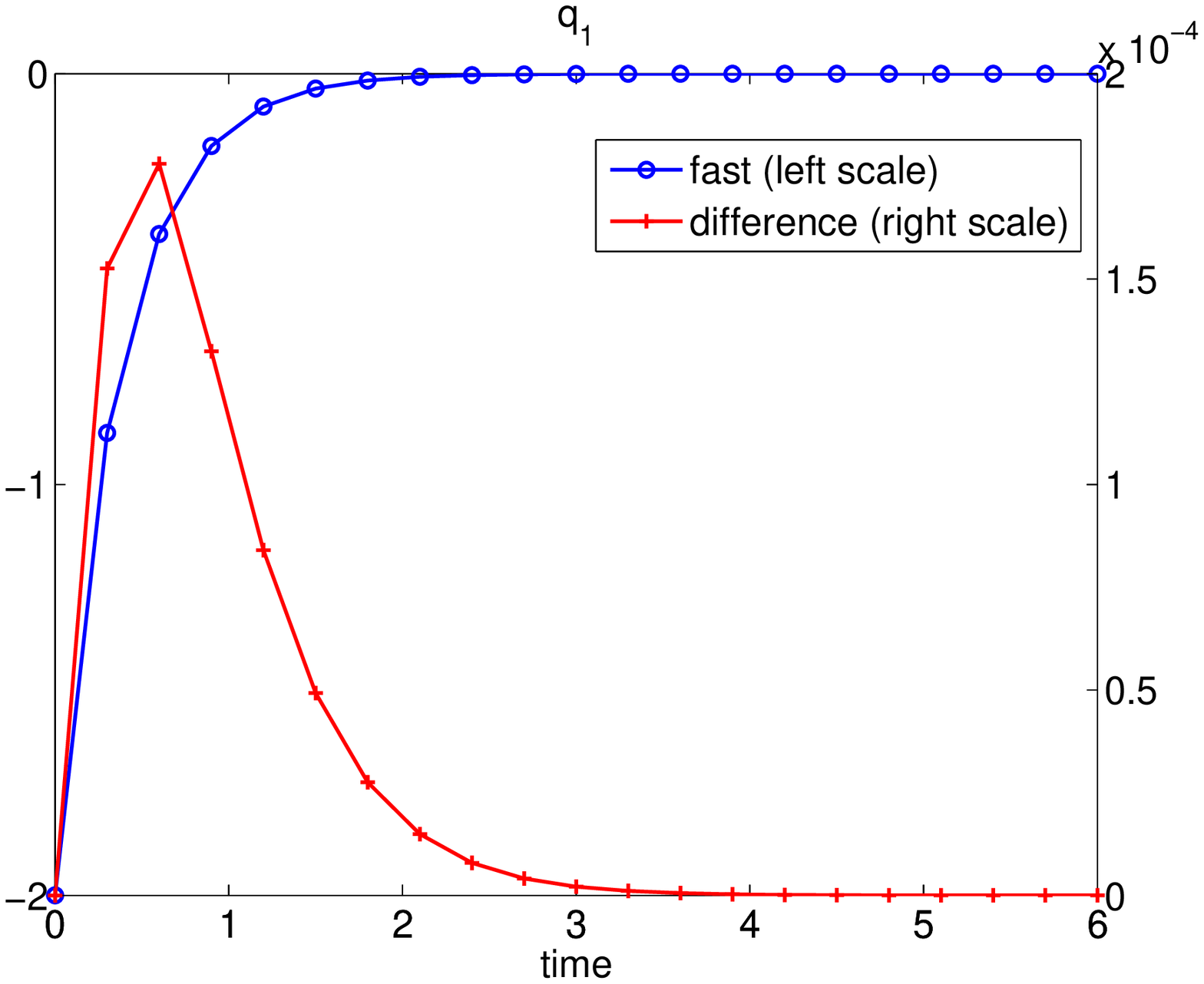}
\includegraphics[width=2.86in]{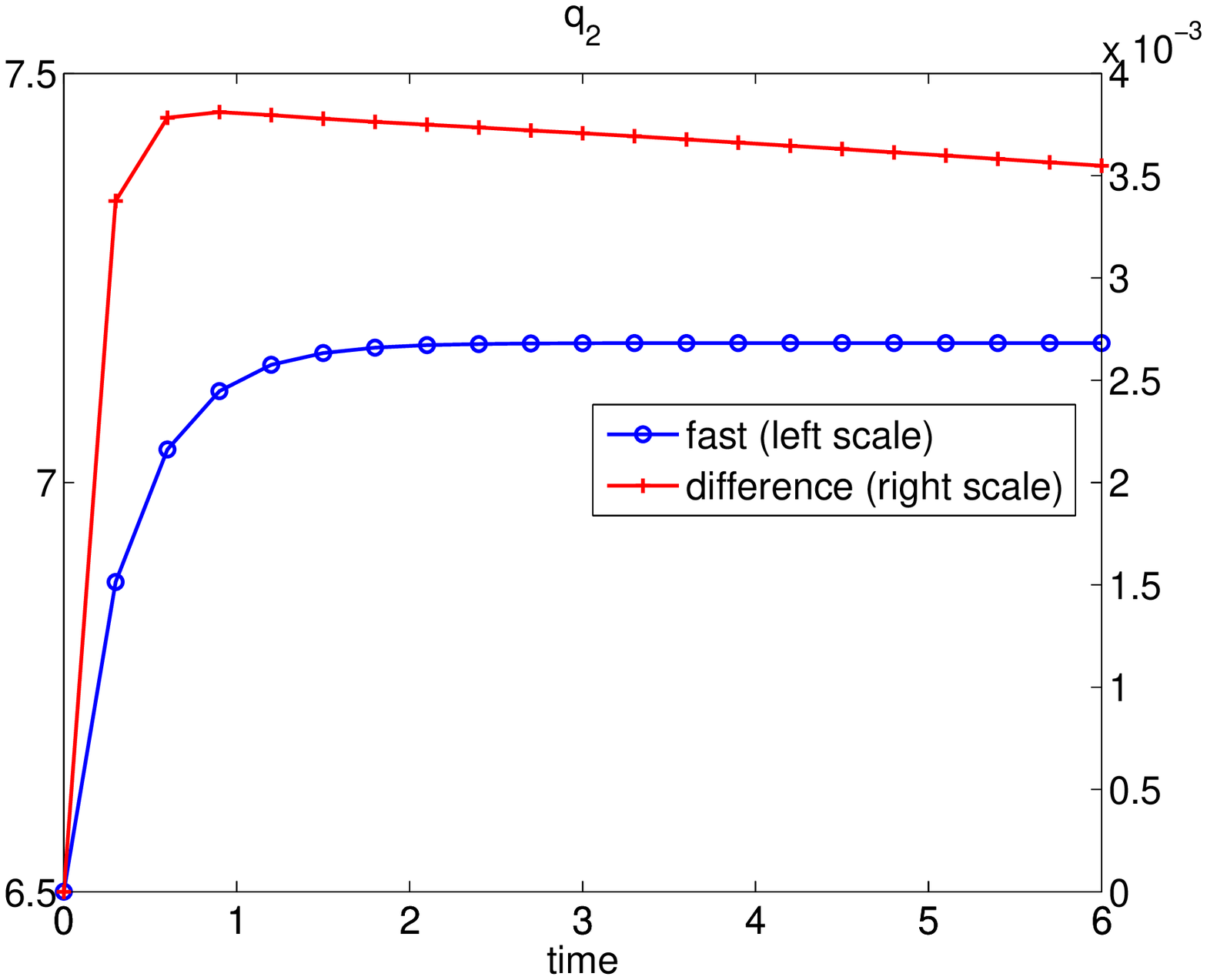}
\end{center}
\caption{Hard sphere molecules. Time evolution for higher order moments. In each figure, the left scale shows the result by the fast spectral method, and the right scale shows the difference between the fast and the direct method. RK4 with $\Delta t=0.3$ for time discretization. $N=32$ in each velocity dimension. In the fast method, $N=32$ in radial direction and $M=74$ for sphere integration. $R=10$, $L=(3+\sqrt{2})R/4\approx11.04$.}
\label{figure4}
\end{figure}

We next consider the same example as in the previous subsection, but for hard sphere molecules. That is, the collision kernel (\ref{VHS}) is assumed to be
\begin{equation}
B=\frac{1}{4\pi}|v-v_*|.
\end{equation}
In this case, there is  no exact formula for either the distribution function or its higher order moments. Therefore, we use the direct spectral method as a reference solution and compare it with the new fast method. The results are plotted in Figure \ref{figure4}, from which we again observe roughly three digits of accuracy for moments. 

\subsection{Angularly dependent collision kernel}

\begin{figure}[htp]
\begin{center}
\includegraphics[width=2.86in]{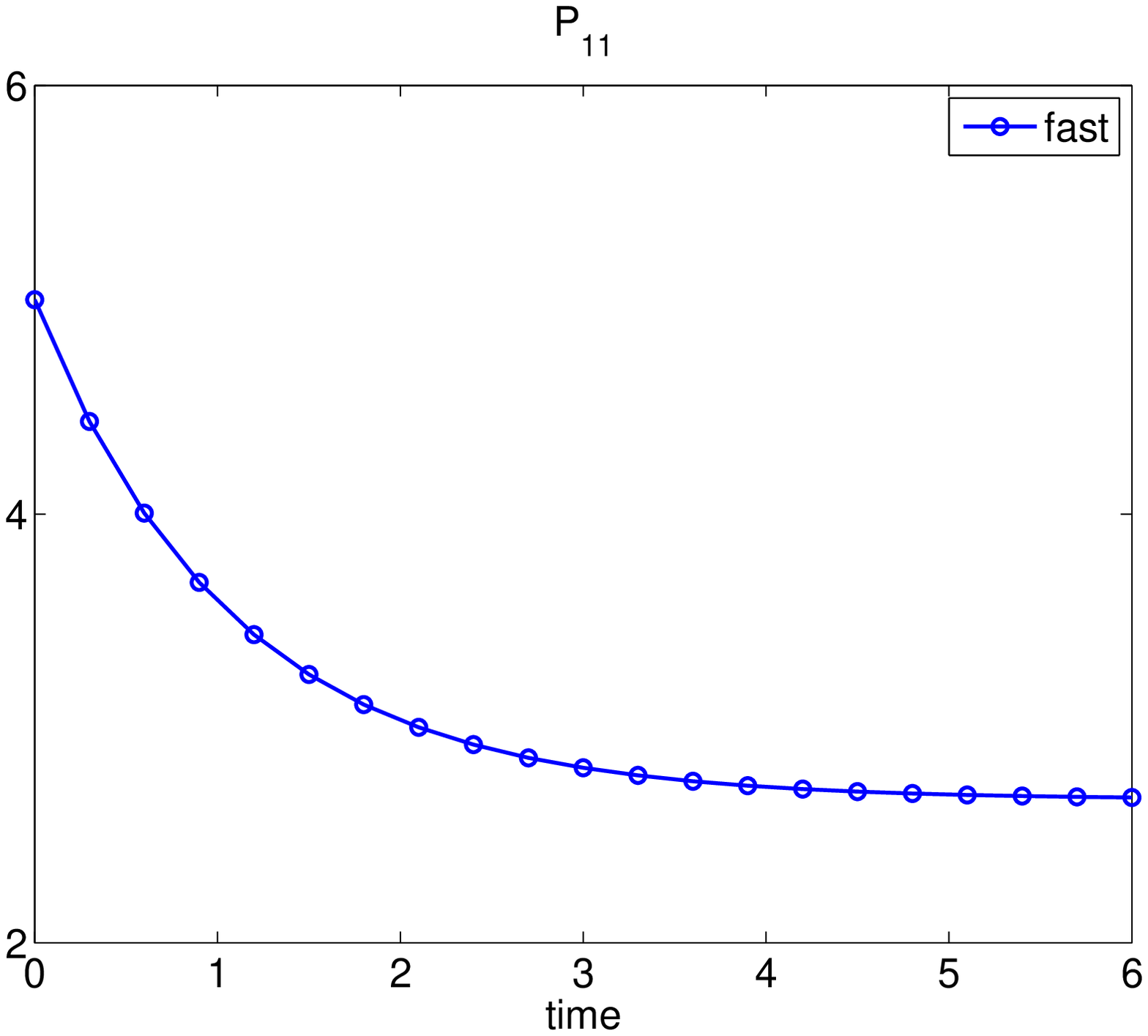}
\includegraphics[width=2.86in]{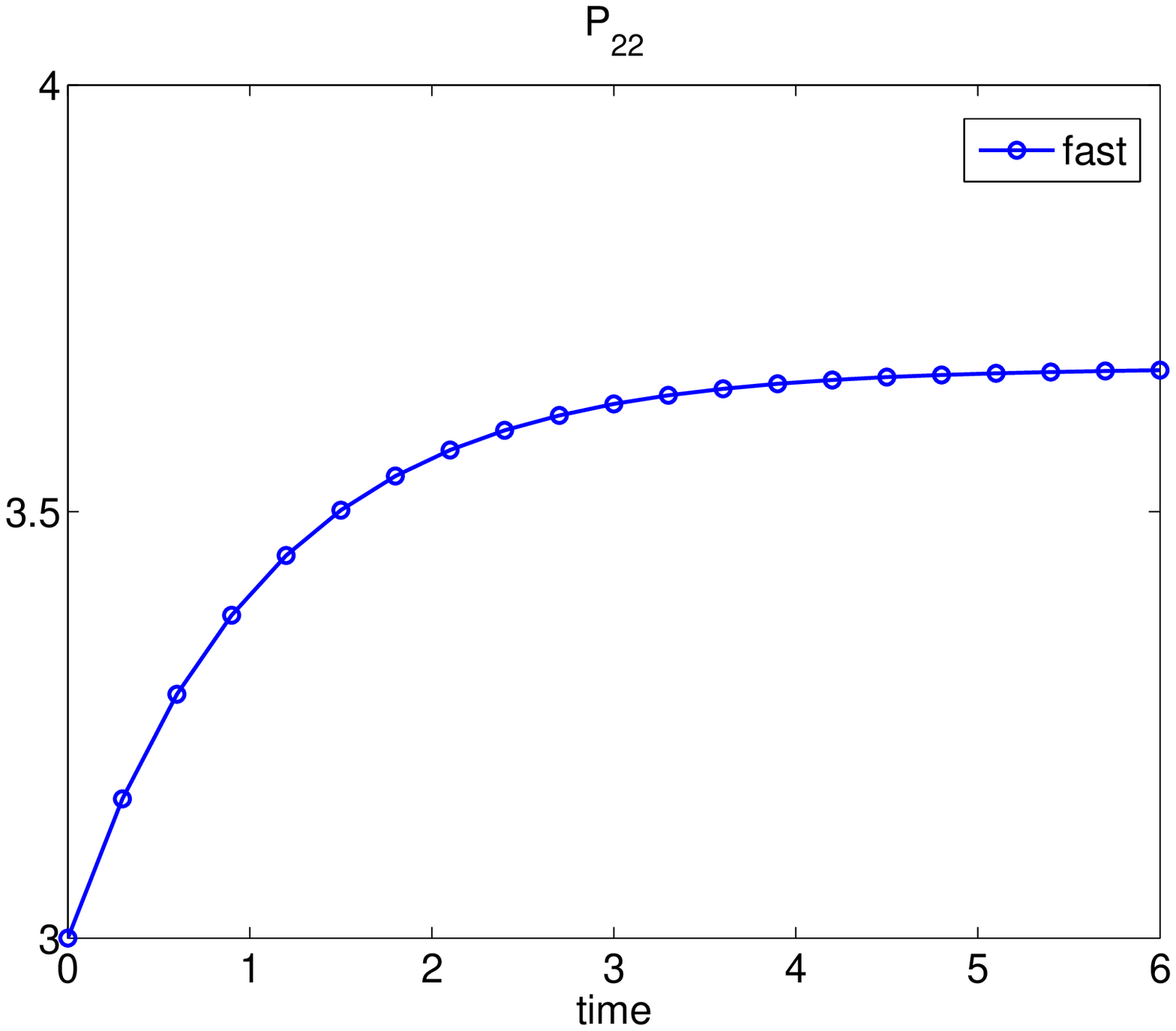}
\includegraphics[width=2.86in]{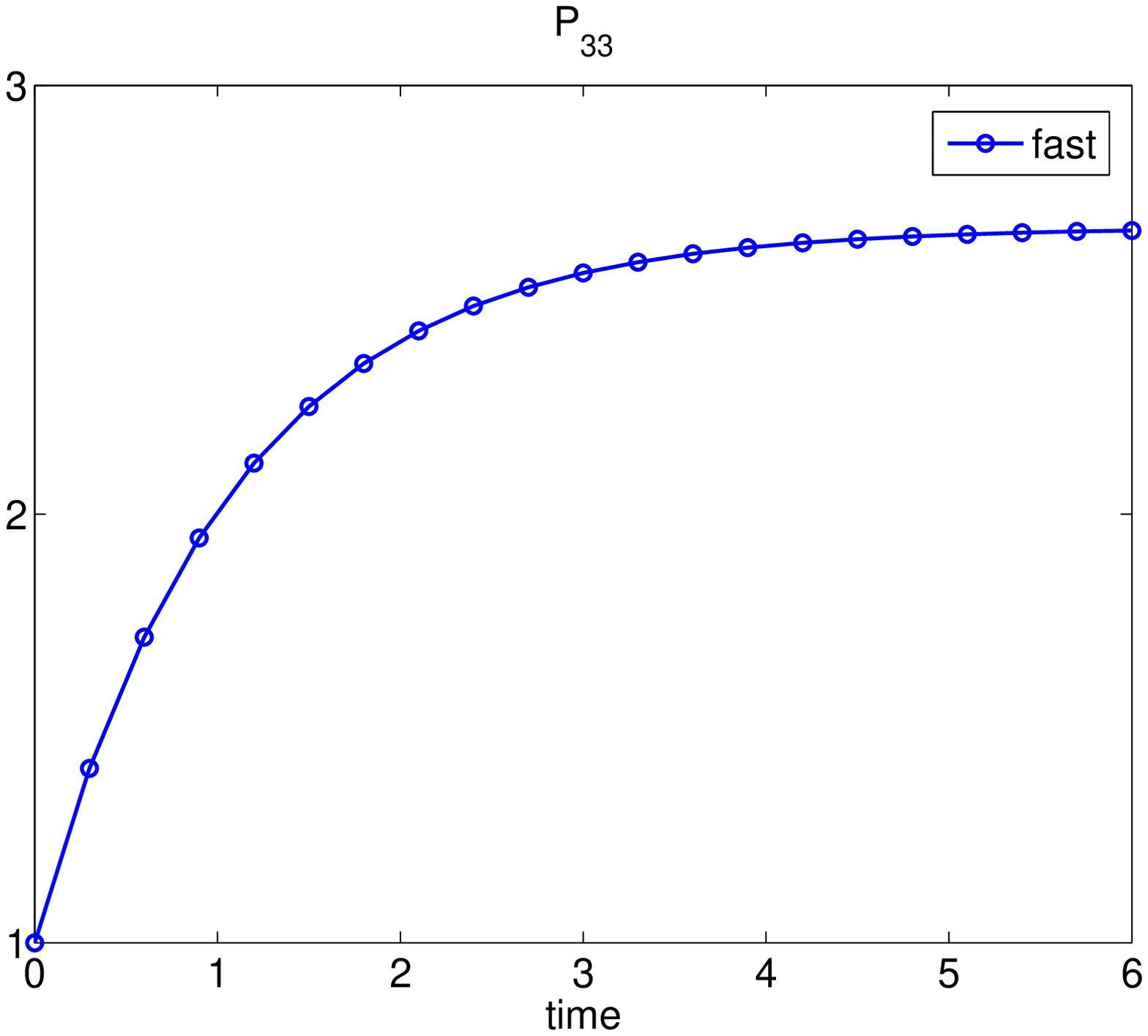}
\includegraphics[width=2.86in]{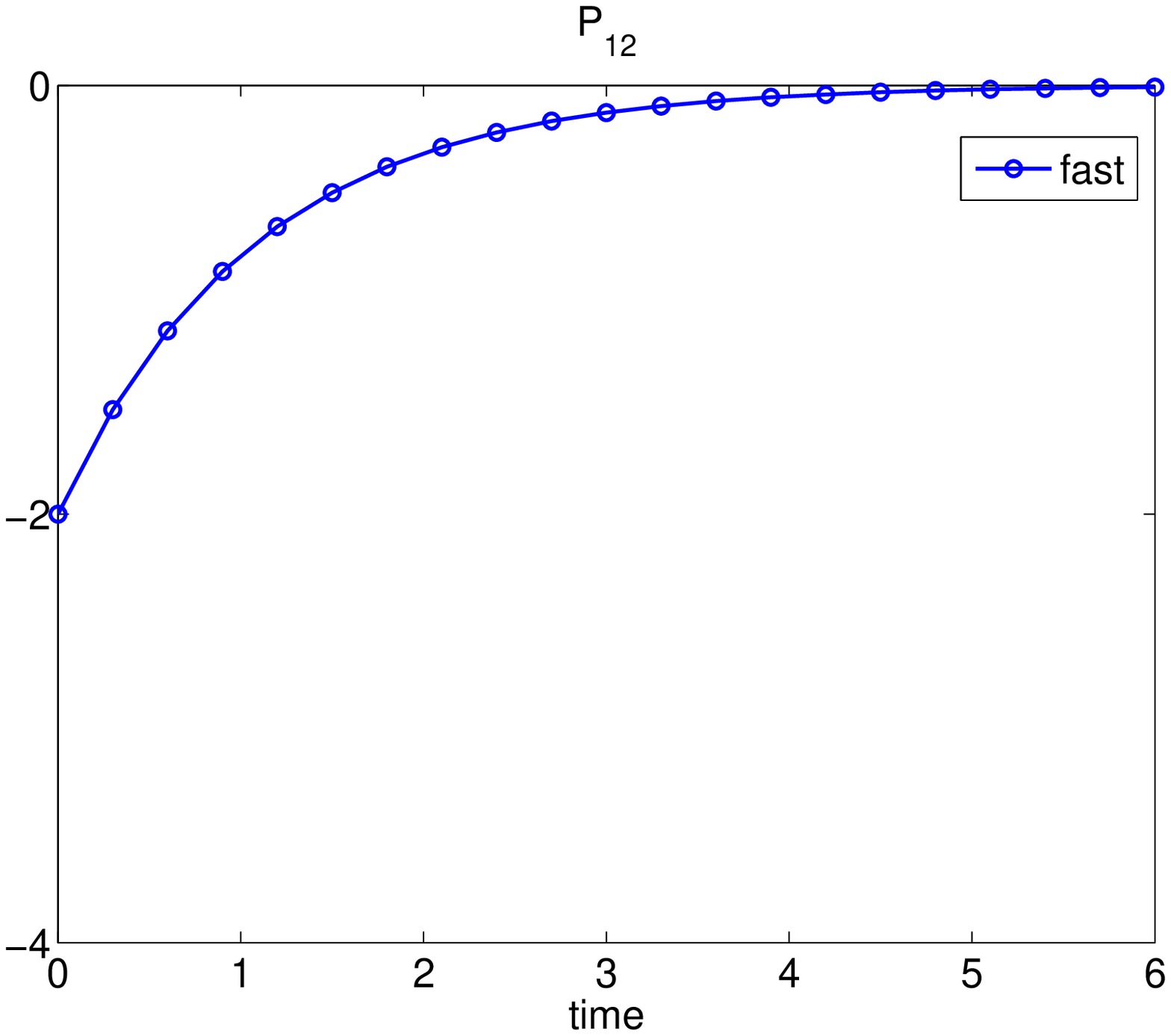}
\includegraphics[width=2.86in]{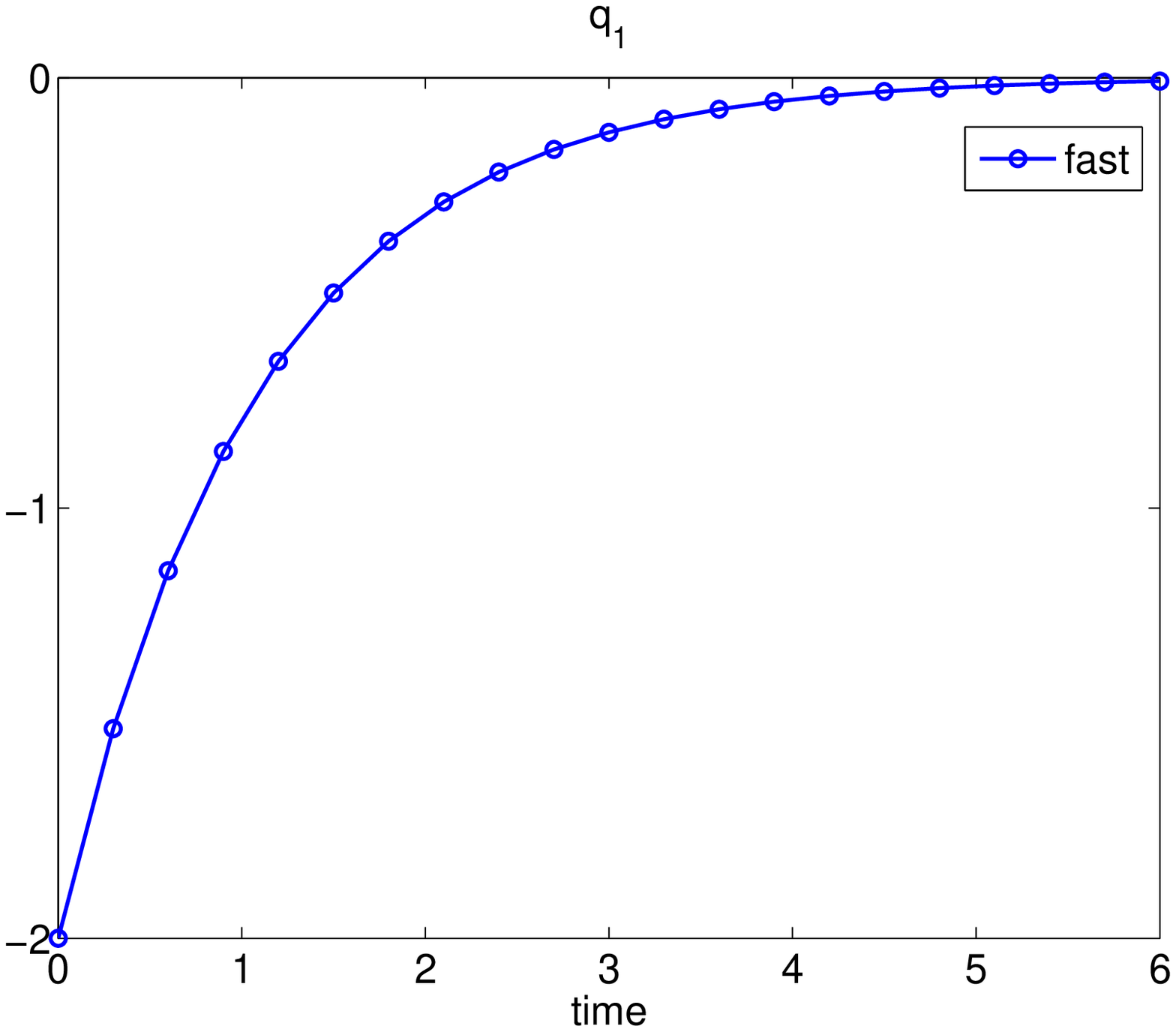}
\includegraphics[width=2.86in]{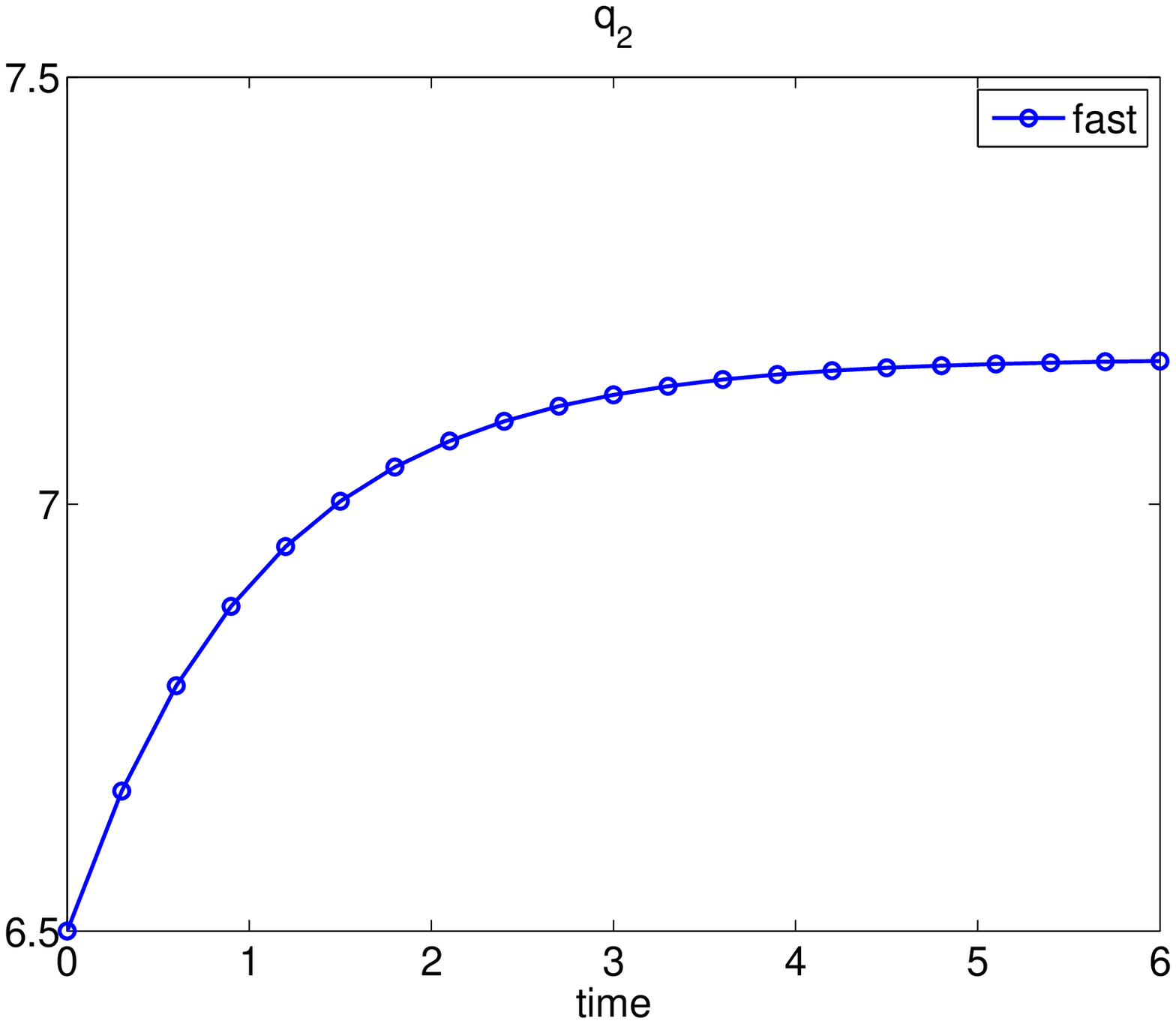}
\end{center}
\caption{Argon molecules. Time evolution for higher order moments computed by the fast spectral method. RK4 with $\Delta t=0.3$ for time discretization. $N=32$ in each velocity dimension, $N=32$ in radial direction and $M=74$ for sphere integration. $R=10$, $L=(3+\sqrt{2})R/4\approx11.04$.}
\label{figure5}
\end{figure}

Our final numerical test involves the variable soft sphere (VSS) model [18], which is widely used in DSMC calculations.  The model has a collision kernel with both velocity and angular dependence:
\begin{equation}
B=b_{\gamma,\eta}|v-v_*|^{\gamma}(1+\cos \theta)^{\eta},
\end{equation}
where $b_{\gamma,\eta}$ is a positive constant and $\cos \theta$ is given in \eqref{CK}.  Setting $\gamma=0.38$, $\eta=0.4$, and $b_{\gamma,\eta} = 1/(4\pi)$, we perform the same test as in Section 4.2 using the same set of discretization parameters.%
\footnote{The choice of $\gamma$ and $\eta$ corresponds to argon gas \cite{Bird} while the choice of $b_{\gamma,\eta}$, which has no effect on the efficiency of the algorithm, is simply a matter of convenience.  In general, these values are composite parameters \cite{WA15} that are used to tune the kernel in order to reproduce experimentally measured values for viscosity and diffusion.  A careful comparison of our method with DSMC for various benchmark examples will be the subject of future work (for which all physical parameters and spatial discretization will need to be included).}
 In Figure \ref{figure5}, we plot the results for the fast method.  Similar results for the direct method are omitted because the time it takes to precompute the weights $G(l,m)$ for this model is prohibitive.  For the fast method, this step takes only a few hours.


\section{Conclusion}
 
A simple, fast spectral method for the Boltzmann collision operator has been proposed in this paper. The method is designed to accelerate the direct method and to relieve the memory bottleneck in its precomputation stage. Through a series of examples, we have demonstrated that the proposed method can be orders of magnitude faster than the direct method while maintaining a comparable level of accuracy.  Furthermore, unlike existing fast spectral methods that can treat only hard sphere molecules, the proposed method is applicable to general collision kernels with both velocity and angular dependence. Ongoing work includes a more careful analysis of spherical quadratures errors and the development of adaptive quadratures to further improve the method.

\section*{Acknowledgements}
J. Hu would like to thank Prof.~Alina Alexeenko for helpful discussion on various collision kernels of practical interest. 
Support from the Institute of Computational Engineering and Sciences (ICES) at the University of Texas Austin is gratefully acknowledged.

\bibliographystyle{plain}
\bibliography{hu_bibtex}
\end{document}